\documentclass{arxivdummy} 
\usepackage{amsmath}
\usepackage{paralist}
\usepackage[misc]{ifsym}
\usepackage{epsfig} 
\usepackage{epstopdf} 
\usepackage[colorlinks=true]{hyperref}

\hypersetup{urlcolor=blue, citecolor=red}
\allowdisplaybreaks

\def\currentvolume{X}
 \def\currentissue{X}
  \def\currentyear{200X}
   \def\currentmonth{XX}
    \def\ppages{X-XX}

\usepackage{amsfonts}
\usepackage{graphicx}
\usepackage{algorithm, algpseudocode} 
\usepackage{cite}
\usepackage{amssymb,amsfonts}
\usepackage{graphicx,bm,url}
\usepackage{textcomp}
\usepackage{booktabs}
\usepackage[capitalise]{cleveref}

\newcommand{\mc}{\mathcal}
\newcommand{\norm}[1]{\left\lVert#1\right\rVert}

\def\x{x_{j+\frac{1}{2}}}
\def\rI{r_{j+\frac{1}{2}}^1}
\def\sI{s_{j+\frac{1}{2}}^1}
\def\rII{r_{j+\frac{1}{2}}^2}
\def\sII{s_{j+\frac{1}{2}}^2}

\def\ueI{\theta_{j+\frac{1}{2}}^1}
\def\ueII{\theta_{j+\frac{1}{2}}^2}
\def\veI{\eta_{j+\frac{1}{2}}^1}
\def\veII{\eta_{j+\frac{1}{2}}^2}

\allowdisplaybreaks

\newtheorem{theorem}{Theorem}[section]

\newtheorem{lemma}[theorem]{Lemma}

\theoremstyle{definition}

\newtheorem{remark}[theorem]{Remark}

\title[Robust Model Reductions of Piezoelectric Beams]
{Robust Model Reductions for the Boundary Feedback Stabilization of Magnetizable Piezoelectric Beams } 

\author[Ahmet Kaan Aydin, Ahmet \"Ozkan \"Ozer, and Jacob Walterman]{}


\thanks{This work was supported by the National Science Foundation of USA under Cooperative Agreement No. 1849213. Additionally, the second author acknowledges the support of the Fulbright U.S. Scholar Program under the 2024-2025 Research Award to conduct research in France.}

\thanks{$^*$Corresponding author: Ahmet Kaan Aydin}

\begin{document}
\maketitle

\centerline{\scshape Ahmet Kaan Aydin$^{{\href{mailto:aaydin1@umbc.edu}{\textrm{\Letter}}}*1}$, Ahmet \"Ozkan \"Ozer$^{{\href{mailto:ozkan.ozer@wku.edu}{\textrm{\Letter}}}2}$ \and Jacob Walterman$^{{\href{mailto:jacob.walterman061@topper.wku.edu}{\textrm{\Letter}}}2}$}

\medskip

{\footnotesize
 \centerline{$^1$Department of Mathematics, University of Maryland, Baltimore County, Baltimore, MD 21250, USA}
}

\medskip

{\footnotesize
 \centerline{$^2$Department of Mathematics, Western Kentucky University, Bowling Green, KY 42101, USA}
} 

\bigskip

\begin{abstract}
Magnetizable piezoelectric beams exhibit strong couplings between mechanical, electric, and magnetic fields, significantly affecting their high-frequency vibrational behavior. Ensuring exponential stability under boundary feedback controllers is challenging due to the uneven distribution of high-frequency eigenvalues in standard Finite Difference models. While numerical filtering can mitigate instability as the discretization parameter tends to zero, its reliance on explicit spectral computations is computationally demanding. This work introduces two novel model reduction techniques for stabilizing magnetizable piezoelectric beams. First, a Finite Element discretization using linear splines is developed, improving numerical stability over standard Finite Differences. However, this method still requires numerical filtering to eliminate spurious high-frequency modes, necessitating full spectral decomposition. Numerical investigations further reveal a direct dependence of the optimal filtering threshold on feedback amplifiers.
To overcome these limitations, an alternative order-reduction Finite Difference scheme is proposed, eliminating the need for numerical filtering. Using a Lyapunov-based framework, we establish exponential stability with decay rates independent of the discretization parameter. The reduced model also exhibits exponential error decay and uniform energy convergence to the original system.
Numerical simulations validate the effectiveness of the proposed methods, and we construct an algorithm for separating eigenpairs for the proper application of the numerical filtering. Comparative spectral analyses and energy decay results confirm the superior stability and efficiency of the proposed approach, providing a robust framework for model reduction in coupled partial differential equation systems.

\end{abstract}

\section{Introduction}

Piezoelectric materials, such as PZT, barium titanate, or lead titanate, are highly versatile smart materials capable of generating electric displacement in direct response to mechanical stress~\cite{smith2005smart}. These materials are critical components in designing actuators that require electrical input, whether in the form of current, charge, or voltage~\cite{ozer2021stabilization,ozer2018potential,yang2009fully}.  The drive frequency is central to their operation, dictating the vibration speed or state changes in piezoelectric beams. Precise control over vibrational modes is achievable by employing periodic or arbitrary signals, offering high resolution and bandwidth. Piezoelectric actuators can achieve positioning resolutions down to tens of nanometers, making them essential in various applications. In addition, piezoelectric materials are extensively used as sensors and energy harvesters, highlighting their multifaceted utility \cite{baur2014advances}.

While electrostatic or quasi-static approximations based on Maxwell's equations are sufficient for describing low-frequency vibrations in many non-magnetizable piezoelectric applications, magnetic effects are often ignored in these approximations. However, magnetic effects can be significant for piezoelectric acoustic wave devices, making existing models inadequate for capturing the vibrational dynamics. Therefore, there is a pressing need for more accurate models that consider the complex interaction between electromagnetic and mechanical interactions. Researchers have stressed the importance of including the complete set of Maxwell's equations as electromagnetic waves come into play. Unlike electrostatic models, the coupling of Maxwell's equations with equations governing mechanical vibrations leads to a fully dynamic theory, often called piezo-electro-magnetism. Researchers have advocated this approach to gain a comprehensive understanding of the behavior of magnetizable piezoelectric materials \cite{darinskii2008role,tallarico2017propagation,yang2009fully,smith2005smart, morris2014modeling,ozer2021stabilization}.

Consider a magnetizable piezoelectric beam clamped on one side and free to oscillate on the other, with a controller regulating tip velocity and tip current accumulation at the electrodes. This beam, length $L$, primarily undergoes longitudinal oscillations, with transverse oscillations considered negligible.

We define positive-definite matrices $\bm C_1$, $\bm C_2$, and $\bm C_3$ as follows \[\bm C_1=\begin{bmatrix}
    \rho &0  \\
  0& \mu        \\
\end{bmatrix} ,\quad \bm C_2=\begin{bmatrix}
  \alpha       & -\gamma \beta   \\
   -\gamma\beta       & \beta  \\
\end{bmatrix},\quad  \bm C_3=\begin{bmatrix}
    k_1 &0  \\
  0& k_2        \\
\end{bmatrix},\]
where $k_1$ and $k_2$ denote positive state feedback gains, and $\rho$, $\alpha$ (and $\alpha_1$), $\beta$, $\gamma$, and $\mu$ represent positive piezoelectric material-specific constants. These constants encompass mass density per unit volume, elastic and piezoelectric stiffness, beam coefficient of impermeability, piezoelectric constant, and magnetic permeability, with $\alpha_1=\alpha-\gamma^2\beta>0.$

Denoting dots $``\cdot"$ as derivatives concerning the time variable $t$, and $v(x,t)$ and $p(x,t)$ as the longitudinal oscillations of the centerline of the beam and the total charge accumulated at the electrodes of the beam, respectively, the equations of motion form a system of strongly-coupled equations
 \begin{equation}\label{eq1}
\begin{cases}
  \bm C_1  \begin{bmatrix}
    \ddot v  \\
    \ddot p       \\
  \end{bmatrix} + \bm C_2 \begin{bmatrix}
  v_{xx}\\
   p_{xx} \\
  \end{bmatrix}=  0,&\\
\begin{bmatrix}
  v\\
  p\\
\end{bmatrix} (0,t)= 0,\quad\left(\bm C_2  \begin{bmatrix}
  v_{x}  \\
   p_{x}       \\
\end{bmatrix} +\bm C_3\begin{bmatrix}
    \dot v \\
    \dot p    \\
\end{bmatrix}\right)(L,t) = 0, &t\in\mathbb{R}^+\\
  \left[v,p,\dot v, \dot p\right]\left(x,0\right)=\left[v_0,p_0,v_1,p_1\right]\left(x\right), & x\in\left[0,L\right].
\end{cases}
\end{equation}
The natural energy of the solutions  of~\eqref{eq1} is defined as
\begin{equation}\label{eq4}
  E(t)=\frac{1}{2}\int^L_0 \left\{ \bm C_1 \begin{bmatrix}
    \dot v  \\
     \dot p       \\
    \end{bmatrix} \cdot   \begin{bmatrix}
    \dot v  \\
     \dot p       \\
    \end{bmatrix}+\bm C_2 \begin{bmatrix}
    v_x  \\
     p_x       \\
    \end{bmatrix} \cdot   \begin{bmatrix}
    v_x  \\
     p_x       \\
    \end{bmatrix} \right\} ~dx,\quad t\ge 0.
\end{equation}

\subsection{Related Literature on the PDE Model}

The PDE system described by~\eqref{eq1} exhibits fundamental limitations regarding exponential stability and exact controllability when only a single state feedback controller is utilized. Specifically, instability emerges when either $k_1=0$ and $k_2\ne 0$, or $k_1\ne 0$ and $k_2=0$. Additionally, certain combinations of material parameters can compromise approximate controllability and strong stability, particularly in scenarios involving high-frequency electromagnetic vibrations~\cite{morris2014modeling}. However, subsequent research by \cite{ozer2015further} indicates that exponential stability can still be attained within a more regularized state space, provided the material parameters are restricted to a narrower range. Crucially, employing both feedback controllers simultaneously ($k_1\ne 0$ and $k_2\ne 0$) ensures exponential stability without imposing restrictive conditions on the material parameters.

The exponential stability of system \eqref{eq1} with both controllers active was first rigorously established by \cite{ramos2019equivalence}, using a decomposition argument combined with an observability inequality. Due to the strong coupling inherent in the system, traditional spectral analysis approaches proved ineffective. Although their method ensured exponential stability, it initially involved suboptimal observation times and did not permit the direct optimization of feedback gains to maximize the decay rate. The observability result was later enhanced in \cite{ozer2022uniform}, where a non-harmonic Fourier series approach provided optimal observation times.

It is important to emphasize that $\sqrt{\frac{\rho}{\alpha_1}}$ and $\sqrt{\frac{\beta}{\mu}}$ (close to the speed of light) represent non-identical wave propagation speeds within the coupled wave equations of system~\eqref{eq1}. Defining the maximum group speed due to strong coupling as
\begin{equation}\label{eta}
\eta:=\max\left(\sqrt{\frac{\rho}{\alpha_1}}+\sqrt{\frac{\mu\gamma^2}{\alpha_1}},\sqrt{\frac{\mu}{\beta}}+\sqrt{\frac{\mu\gamma^2}{\alpha_1}}\right),
\end{equation}
and introducing the Hilbert space $\mathcal H:=(H^1{(0,L)})^2\times (L^2(0,L))^2$, with $(H^1(0,L))=\{z\in H^1(0,L): z(0)=0\}$, a refined analysis was conducted.

Addressing the limitations of prior methodologies, recent research by \cite{ozer2024exponential} introduced a carefully constructed Lyapunov function. This method explicitly determines the maximal decay rate and directly enables the optimization of feedback gains, thereby eliminating dependence on observability inequalities. The principal result is summarized in the following theorem:

\begin{theorem}\label{thm1}\cite[Theorem 4]{ozer2024exponential}
The energy $E(t)$ is dissipative; specifically, $$\frac{dE(t)}{dt}=-k_1\left|v_t(L,t)\right|^2-k_2\left|p_t(L,t)\right|^2\le 0$$ for all $t>0$. Moreover, for any $\epsilon>0$ and initial conditions $(v,p,\dot v,\dot p),(v_0,p_0,v_1,p_1)\in \mathcal{H}$, the energy $E(t)$ decays exponentially, satisfying
\begin{align*}
  &\begin{cases}
    E(t)\leq M E(0)e^{-\sigma t},&\forall t>0,\\
\sigma(\delta(\xi_1,\xi_2,\epsilon))= \delta\left(1-\delta L \eta\right),\quad M(\delta(\xi_1,\xi_2,\epsilon))=\frac{1+\delta L \eta}{1-\delta L\eta}, & \text{where}
  \end{cases}\\
  &\begin{cases}
    \delta(k_1,k_2,\epsilon):=\frac{1}{L}\text{min}\left(\frac{1}{\eta}, f_1(k_1,\epsilon), f_2(k_2,\epsilon)\right),\\
      f_1(k_1,\epsilon):=\frac{2k_1\alpha_1}{\rho\alpha_1+(1+\epsilon)k_1^2},\quad f_2(k_2,\epsilon):=\frac{2k_2\epsilon\alpha_1\beta}{\epsilon\mu\alpha_1\beta+(\epsilon\alpha+\gamma^2\beta)k_2^2}.
\end{cases}
\end{align*}
\end{theorem}

From a practical perspective, determining appropriate intervals for feedback amplifier gains $k_1$ and $k_2$ is essential to achieving the targeted decay rate $-\sigma$ indicated in Theorem~\ref{thm1}. The maximal decay rate is obtained when $\delta = \frac{1}{2\eta L}$, yielding $\sigma_{\text{max}} = \frac{1}{4\eta L}$. In \cite{ozer2024exponential}, a multiplier method coupled with optimization arguments explicitly identifies these intervals, determining the maximal achievable decay rate. This rate provides an upper bound for the optimal decay rate, which may be exceeded depending on specific parameter configurations.

Recent literature has further expanded methods addressing system \eqref{eq1}. A passivity-based approach introduced in \cite{de2023control} and various damping configurations have been investigated. Specifically, distributed damping approaches established exponential stabilization with a single viscous damping element \cite{ramos2018exponential}. Moreover, exponential stability has been demonstrated using two delayed-type damping designs across various boundary and distributed damping configurations, as shown in \cite{feng2022exponential, kong2022equivalence}.

Recent studies have explored fractional damping and thermal effects in nonlinear piezoelectric beam designs, leading to reductions into a single boundary-controlled wave equation known for exponential stability \cite{freitas2022long}. Furthermore, investigations have incorporated long-range boundary and distributed memory terms into system \eqref{eq1}, achieving exponential stability \cite{feng2023stability} and polynomial stability \cite{zhang2023stability}. Moreover, output feedback stabilization and non-collocated observer designs have been emphasized to attain exponential stability \cite{ozer2024boundary}. Finally, recent results on exponential and polynomial stability of piezoelectric beams in transmission-line settings with partial damping designs are provided in \cite{akil2024stability,akil2025advancing}.

\subsection{Related Literature on Model Reductions}

Model reduction plays a crucial role in the simulation and control of PDE models, typically achieved through numerical discretization methods such as Finite Differences or Finite Elements, which convert PDEs into finite-dimensional ordinary differential equations (ODEs). However, these finite-dimensional approximations often introduce unintended high-frequency modes, resulting in the well-known ``spillover effect" \cite{preumont2018vibration}. Additionally, standard discretization methods frequently produce spurious and unobservable vibrational modes, complicating accurate modeling and control \cite{banks2012functional, banks1991exponentially, zuazua2006controllability}.

To overcome these issues, filtering techniques such as \textit{direct Fourier filtering} and \textit{indirect filtering} have been developed \cite{infante1999boundary, tebou2007uniform}. Direct Fourier filtering effectively isolates eigenvalues corresponding to low and high frequencies, removing undesirable modes while preserving convergence. This method has been successfully applied to diverse PDE models, including wave equations, Euler-Bernoulli beams, Rayleigh beams, and multi-layer beam systems, using decomposition and Lyapunov-based approaches \cite{infante1999boundary, el2013boundary, ozer2024exponential, leon2002boundary, ozer2019uniform, ozer2022robust}.

However, direct Fourier filtering is mainly applied in the literature for single, decoupled, or diagonalizable PDE systems with distinctly separated branches of eigenvalues. Coupled PDE systems, featuring multiple potentially entangled eigenvalue branches, require separate and meticulous filtering of each branch. Differentiating and filtering these branches become notably challenging, especially in control applications involving nonzero feedback amplifiers.

A Finite-Difference-based model reduction specifically addressing boundary observations of system \eqref{eq1} with zero feedback parameters ($k_1=k_2=0$) was introduced in \cite{ozer2022uniform}. This approach employs direct Fourier filtering uniformly across two branches of eigenvalues, ensuring uniform observability as the discretization parameter approaches zero. Special attention is required when filtering high-frequency eigenvalues related to magnetic effects, as indiscriminate filtering can unintentionally eliminate important modes.

An additional limitation of direct Fourier filtering is its inherent ambiguity regarding optimality, which raises concerns about achieving the best possible model reduction \cite{lissy2019optimal}. Recently, advanced Finite-Difference-based model reduction techniques have emerged for PDEs, such as one-dimensional wave equations \cite{liu2020new} and Euler-Bernoulli beam models \cite{aydin2023novel,aydinhaider2024, liu2019novel, liu2021uniformly}. These methods construct reduced-order models using equidistant grid points and averaging operators, completely avoiding numerical filtering. The uniform observability of these reduced models is demonstrated directly through discrete multiplier methods, offering a robust alternative to classical filtering approaches.

\subsection{Our Contributions and Content}

In this paper, we present novel contributions to the field of model reduction for the system \eqref{eq1}. Specifically, we introduce a Finite Elements-based model reduction using linear splines, an approach not previously explored for this class of coupled PDE systems. Our analysis reveals that this model reduction lacks exponential stability as the discretization parameter approaches zero unless proper numerical filtering is applied. However, with appropriate filtering, this method provides a robust model reduction framework, requiring significantly less filtering compared to standard Finite Differences-based reductions. A key challenge in this approach is the necessity of computing the full spectral decomposition of the system for Fourier filtering, which becomes computationally expensive for large-scale discretizations. As a result, one open problem is the development of a more efficient algorithm that filters out spurious modes without requiring explicit eigenvalue computation.

Since direct Fourier filtering must be applied separately to each branch of eigenvalues, we introduce a novel algorithm to efficiently distinguish and separate these branches. This algorithm, detailed in Remark \ref{rmk:filtering}, provides an automated approach for ensuring accurate application of numerical filtering, addressing a major computational challenge in controlled coupled PDE systems.

Motivated by recent order-reduction methodologies in \cite{aydin2023novel, aydinhaider2024,liu2019novel,liu2021uniformly}, we further propose a Finite Differences-based model reduction approach that eliminates the need for numerical Fourier filtering altogether. This advancement is particularly valuable for systems with multiple eigenvalue branches, such as \eqref{eq1}, where filtering complexity and ambiguity—especially regarding high-frequency magnetic effects—pose significant challenges. By bypassing the need for filtering, this approach enhances numerical stability, robustness, and computational efficiency.

Finally, we establish exponential stability and energy convergence results for the proposed reduced-order models using a Lyapunov framework. To the best of our knowledge, this work represents the first demonstration of exponential decay for error dynamics and uniform convergence of discretized energy in the $C[0, \infty]$-norm to the original PDE model’s energy for this system. Notably, the decay rate of our reduced-order model (Theorem \ref{mainexp}) precisely matches that of the original PDE model (Theorem \ref{thm1}), previously established in \cite{ozer2024exponential}. Consequently, the optimization strategies developed for the PDE model in \cite{ozer2024exponential} directly apply to the reduced model, significantly enhancing practical applicability. These results offer deeper insights into the stability and accuracy of reduced models for strongly coupled wave-equation systems, such as those investigated in \cite{ozer2017modeling, ozer2018potential, ozer2021stabilization}.

A key emphasis of this paper is on numerical implementation and validation. In particular, we develop and implement an algorithmic framework for model reduction and filtering, ensuring computational feasibility and efficiency for large-scale discretizations. The numerical section not only validates the theoretical predictions but also presents novel insights into the impact of filtering, feedback amplifiers, and eigenvalue separation on stability performance. Additionally, our findings indicate a strong interplay between the feedback amplifiers and the optimal filtering parameter. Simulations reveal that different feedback gain selections alter the necessary amount of filtering required to achieve the fastest decay rate, highlighting an intricate relationship between control design and spectral properties of the reduced model. However, a precise theoretical characterization of this dependency remains an open research direction.

The paper is organized as follows. Section \ref{FEMSplines} introduces the Finite Element discretization using linear splines and examines its instability in the absence of numerical filtering. Section \ref{ModelRed} presents the Finite Differences-based model reduction, which eliminates numerical Fourier filtering and includes theoretical analyses on exponential stability, error dynamics, and energy convergence while also investigating the impact of feedback amplifiers on the optimal filtering threshold. Section \ref{Simu} is dedicated to numerical simulations and computational algorithms, introducing a new algorithm for separating the branches of eigenvalues on the left-half plane, and showcasing the dependency of the optimal filtering threshold on feedback amplifiers for Finite Element discretization. In addition to validating theoretical results, comparative spectral analyses, energy decay rates, and eigenvalue distributions further illustrate the superior stability of novel Finite Differences based model over Finite Elements. Section \ref{conc} concludes with a summary of findings and future research directions.

\section{Model Reduction by the Finite Elements Method (FEM) with Linear Splines}\label{FEMSplines}
To derive a reduced model corresponding to \eqref{eq1}, we employ the Finite Elements Method (FEM). Let $N \in \mathbb{N}$ be given, and define the mesh size as $h := \frac{L}{N+1}$. Consider a uniform discretization of the interval
 $$0=x_0<x_1<...<x_i=i*h<\ldots <x_{N+1}=L.$$
First, multiply both sides of the equation in \eqref{eq1} by a continuously differentiable test function $\bm\phi(x) = \begin{bmatrix} \phi^1(x) \ \phi^2(x) \end{bmatrix} \in (C_0^\infty[0,L])^2$, and integrate over $[0, L]$. This yields the following weak formulation
 \begin{equation}\label{weak}
 \bm C_1  \int_0^L\begin{bmatrix}
  \ddot v  \\
  \ddot p       \\
\end{bmatrix} \cdot \begin{bmatrix}
\phi^1\\
  \phi^2       \\
\end{bmatrix}dx- \bm C_2  \int_0^L\begin{bmatrix}
v_{x}  \\
 p_{x}       \\
\end{bmatrix} \cdot \begin{bmatrix}
\phi^1_{x}  \\
\phi^2_x     \\
\end{bmatrix}dx+ \bm C_3  \begin{bmatrix}
  \dot v \\
  \dot p      \\
\end{bmatrix} \cdot  \begin{bmatrix}
  \phi^1  \\
  \phi^2       \\
\end{bmatrix}(L)= \bm 0.
\end{equation}
For each $j = 1, 2$ and at each node ${x_i}{i=1}^N$, the linear basis splines are defined as follows
\begin{equation*}
  \phi_i^j(x)= \begin{cases}
    \frac{x-x_i}{h},& x_{i-1}<x<x_{i}\\
\frac{x_{i+1}-x}{h},& x_{i}<x<x_{i+1}\\
0,& {\textrm otherwise}
  \end{cases},\quad \phi^j_{N+1}(x)= \begin{cases}
  \frac{x-x_N}{h},& x_{N}<x<x_{N+1}\\
0,& {\textrm otherwise.}
  \end{cases}
\end{equation*}

Let $(v_i, p_i) = (v_i, p_i)(t) \approx (v, p)(x_i, t)$ denote the approximation of the solution $(v, p)(x, t)$ of \eqref{eq1} at the nodes $x_i$ for $i = 0, 1, \ldots, N, N+1$. Define the vectors $\bm{v} = [v_1, v_2, \ldots, v_{N+1}]^T$ and $\bm{p} = [p_1, p_2, \ldots, p_{N+1}]^T$. Using central differences, approximate $v_{xx}(x_i) \approx (-\bm{A}_h \bm{v})_i$ and $p{xx}(x_i) \approx (-\bm{A}_h \bm{p})_i$. Define the $(N+1) \times (N+1)$ matrices $\bm{A}_h$, $\bm{M}$, and $\bm{B}$ as
\begin{align}
	\bm A_h&:=\frac{1}{h^2}\begin{bmatrix}
		2&-1&0&\dots&\dots&\dots&0\\
		-1&2&-1&0&\dots&\dots&0\\
		&\ddots&\ddots&\ddots&\ddots&\ddots&\\
		0&\dots&\dots&0&-1&2&-1\\
		0&\dots&\dots&\dots&0&-1&1\\
	\end{bmatrix},\nonumber \\ \label{eq:massFEM}
  \bm M&:=\frac{1}{6}\begin{bmatrix}
		4&1&0&\dots&\dots&\dots&0\\
		1&4&1&0&\dots&\dots&0\\
		&\ddots&\ddots&\ddots&\ddots&\ddots&\\
		0&\dots&\dots&0&1&4&1\\
		0&\dots&\dots&\dots&0&1&2\\
	\end{bmatrix}, \quad \bm B:=\begin{bmatrix}
		0&0&\dots&0\\
		0&0&\dots&0\\
		\vdots&\ddots&\ddots&\vdots&\\
		0&\dots&0&1
	\end{bmatrix}.
	\end{align}
The weak formulation \eqref{weak} is thereby reduced to
	\begin{equation}
	\label{FEM-pseudo}
	~\left\{	\begin{array}{ll}
(\bm C_1 \otimes \bm M) \begin{bmatrix}
     {\ddot{\bm{v}}}   \\
   {\ddot{\bm{p}}}        \\
\end{bmatrix} +\left(\bm C_2  \otimes \bm A_h\right) \begin{bmatrix}
     {{\bm{v}}}   \\
   {{\bm{p}}}        \\
\end{bmatrix}	 + (\bm C_3 \otimes \bm B) \begin{bmatrix}
     {\dot {\bm{v}}}   \\
   {\dot {\bm{p}}}        \\
\end{bmatrix} =\bm 0,  & t\in \mathbb{R}^+,\\
(\bm v, \dot{\bm{v}},\bm p,\dot{\bm p})_i(0)=(v^0,v^1,p^0, p^1)(x_i),  & i=0,...,N+1.
	\end{array}\right.
\end{equation}
Here, $\otimes$ represents the Kronecker product. Using the mixed-product property of the Kronecker product, the system \eqref{FEM-pseudo} can be further simplified as
\begin{eqnarray}
	\label{FEM}
 \begin{bmatrix}
     {\ddot{\bm{v}}}   \\
   {\ddot{\bm{p}}}        \\
\end{bmatrix} +\left(\bm C_1^{-1}\bm C_2\right)  \otimes \left(\bm M^{-1}\bm A_h\right) \begin{bmatrix}
     {{\bm{v}}}   \\
   {{\bm{p}}}        \\
\end{bmatrix}
	 + \left(\bm C_1^{-1}\bm C_3\right) \otimes \left(\bm M^{-1} \bm B\right) \begin{bmatrix}
     {\dot {\bm{v}}}   \\
   {\dot {\bm{p}}}        \\
\end{bmatrix} =\bm 0, t\in \mathbb{R}^+.
\end{eqnarray}

Let $\bm Z_h = \frac{1}{h}{\rm tridiag}(0, -1, 1)$ denote an $(N+1)\times (N+1)$ tridiagonal matrix. The discretized energy corresponding to \eqref{FEM}, which serves as the discretized counterpart of \eqref{eq4}, is expressed as
\begin{align*}
E_h^{\rm FEM}(t):=&\frac{h}{2}\left\{(\bm C_1 \otimes \bm M)  \begin{bmatrix}
   {\dot {\bm v}}\\
{\dot {\bm p}}     \\
\end{bmatrix} \right\}\cdot \begin{bmatrix}
   {\dot {\bm v}} \\
{\dot {\bm p}}     \\
\end{bmatrix}+\frac{h}{2}\left\{(\bm C_2 \otimes\bm  Z_h)  \begin{bmatrix}
   { {\bm v}} \\
{ {\bm p}}     \\
\end{bmatrix} \right\}\cdot \left\{(\bm I \otimes\bm  Z_h)  \begin{bmatrix}
   { {\bm v}} \\
{ {\bm p}}     \\
\end{bmatrix} \right\}\\
& +\frac{h}{12}\left(\bm C_1 \begin{bmatrix}
   {\dot {v}}_{N+1} \\
{\dot {p}}_{N+1}     \\
\end{bmatrix}\cdot \begin{bmatrix}
   {\dot {v}}_{N+1} \\
{\dot {p}}_{N+1}     \\
\end{bmatrix} +6\bm C_2 \begin{bmatrix}
   {{v}}_{N+1} \\
{{p}}_{N+1}     \\
\end{bmatrix}\cdot \begin{bmatrix}
   { {v}}_{N+1} \\
{{p}}_{N+1}     \\
\end{bmatrix}\right).
\end{align*}

\subsection{Spectrum and Fourier Solutions of the Control-Free System}

Let $\bm\Phi(t) = \begin{bmatrix} \bm v & \bm p & \dot{\bm v} & \dot{\bm p} \end{bmatrix}^T.$ The control-free system, \eqref{FEM} with $k_1, k_2 \equiv 0$, i.e., $\bm C_3 \equiv \bm 0$, can be rewritten in first-order form as
\begin{equation}
  \label{disc}
  \dot{\bm\Phi}= {\mc A}^{\rm FEM} \bm\Phi :=
  \begin{bmatrix}
    \bm 0& \bm I\\
    \mathbb{A}& \bm 0
  \end{bmatrix}\bm\Phi, \quad \begin{cases}
    \mathbb{A}= -\left(\bm C_1^{-1}\bm C_2\right)  \otimes \left(\bm M^{-1}\bm A_h\right),\\
    \bm\Phi(0) = \begin{bmatrix}
      \bm v^0 & \bm p^0 & \bm v^1 & \bm p^1
    \end{bmatrix}^T.
  \end{cases}
  \end{equation}
Next, consider the eigenvalue problem associated with \eqref{disc}
\begin{equation}
  \label{eq:EigProbFEM}
{{\mc A}}^{\rm FEM}  \bm \Psi={{\tilde \lambda} }(h)\bm \Psi.
 \end{equation}
The eigenvalues of the Kronecker product of two matrices are the product of the eigenvalues of the individual matrices, and the eigenvectors are the Kronecker product of the eigenvectors of the individual matrices. Furthermore, the eigenpairs of ${{\mc A}}^{\rm FEM}$ can be derived from the eigenpairs of $\mathbb{A}$. Therefore, it suffices to compute the eigenpairs of $\bm C_1^{-1}\bm C_2$ and $\bm M^{-1}\bm A_h$ to solve the eigenvalue problem \eqref{eq:EigProbFEM}.

\begin{theorem}\cite[Lemma 3]{ozer2023revisiting} For $k =1,  \ldots, N+1,$ the eigenvalues, $\lambda_k$ of $\bm M^{-1}\bm A_h$  are given by
  \begin{equation}
    \label{eq:lambda}
  {\lambda}_k =\frac{1}{h^2}\frac{6-6\cos{\left(\frac{(2k-1)\pi h}{2(L-h)}\right)}}{2+\cos{\left(\frac{(2k-1)\pi h}{2(L-h)}\right)}}.
  \end{equation}
  Moreover, $\lambda_k h^2<12$  for all $h>0,$ and $\lambda_{N+1} h^2\rightarrow 12$ as $h\rightarrow 0.$
  \end{theorem}
\noindent Define the following constants  $b_1 := \frac{1}{\gamma\mu}(\frac{\alpha_1}{\zeta_1^2}-\rho)$, $b_2 := \frac{1}{\gamma\mu}(\frac{\alpha_1}{\zeta_2^2}-\rho),$
\begin{align*}
\zeta_1 &:= \frac{1}{\sqrt{2}}\sqrt{\frac{\alpha}{\rho}+\frac{\beta}{\mu}
+\sqrt{\left|\frac{\alpha}{\rho}+\frac{\beta}{\mu}\right|^2-\frac{4\alpha_1}{\mu\rho}}}, \\
 \zeta_2 &:=  \frac{1}{\sqrt{2}}\sqrt{\frac{\alpha}{\rho}+\frac{\beta}{\mu}
-\sqrt{\left|\frac{\alpha}{\rho}+\frac{\beta}{\mu}\right|^2-\frac{4\alpha_1}{\mu\rho}}},
 \end{align*}
with $\zeta_1, \zeta_2 > 0$, $b_1, b_2 \neq 0$, $b_1 \neq b_2$, and $b_1b_2 = -\frac{\rho}{\mu}$. Additionally, $b_1$ and $b_2$ satisfy the quadratic equation
\begin{equation} \label{b-quad} b^2 - \left(\frac{\alpha}{\gamma\beta} - \frac{\rho}{\gamma\mu}\right)b - \frac{\rho}{\mu} = 0.
\end{equation}

The eigenvalues of the matrix $\left(\bm C_1^{-1}\bm C_2\right)$ are easily calculated as ${\zeta_1^2, \zeta_2^2}.$ Therefore, the spectrum of $\tilde{\mc A}^{\rm FEM}$ can be constructed as shown in the following theorem and summarized in Table \ref{tab:eigenvalues_eigenvectors}.
\begin{theorem} 
Assume $\bm  C_3\equiv \bm 0$ in \eqref{FEM}. The  matrix $\tilde{\mc A}^{\rm FEM}$  has two branches of eigenvalues
 \begin{equation}\label{EVs}
 \left\{\tilde \lambda_{1k}^\mp(h)=\mp i \zeta_1\sqrt{{\lambda}_k },\quad
  \tilde \lambda_{2k}^\mp(h)=\mp i \zeta_2\sqrt{{\lambda}_k },
 \right.
 \end{equation}
 and since  $\tilde \lambda_{1k}^-=-\tilde\lambda_{1k}^+, ~~\tilde\lambda_{2k}^-=-\tilde \lambda_{2k}^+,$
the corresponding eigenvectors are
\begin{eqnarray*}
{\bm\Psi}^+_{1k}(h)=\left( \begin{array}{c}
 \frac{1}{\tilde\lambda_{1k}^+} {\bm\psi}_k\\
 \frac{b_1}{\tilde\lambda_{1k}^+}{\bm\psi}_k \\
 {\bm\psi}_k\\
b_1 {\bm\psi}_k
 \end{array} \right), ~
 {\bm\Psi}^{-}_{1k}(h)=\left( \begin{array}{c}
 \frac{1}{\tilde\lambda_{1k}^+}{\bm\psi}_k \\
 \frac{b_1}{\tilde \lambda_{1k}^+} {\bm\psi}_k\\
 -{\bm\psi}_k\\
-b_1{\bm\psi}_k
 \end{array} \right),\\
~ {\bm\Psi}^+_{2k}(h)=\left( \begin{array}{c}
 \frac{1}{\tilde\lambda_{2k}^+}{\bm\psi}_k\\
 \frac{b_1}{\tilde \lambda_{2k}^+} {\bm\psi}_k\\
 {\bm\psi}_k\\
b_2 {\bm\psi}_k
 \end{array} \right), ~
 {\bm\Psi}^{-}_{2k}(h)=\left( \begin{array}{c}
 \frac{1}{\tilde\lambda_{2j}^+} {\bm\psi}_k\\
 \frac{b_2}{\tilde \lambda_{2j}^+} {\bm\psi}_k\\
-{\bm\psi}_k\\
-b_2 {\bm\psi}_k
 \end{array} \right).\end{eqnarray*}
  where $  { \psi}_{k,j}=sin\left(\frac{(2k-1)j\pi}{2N}\right),$ for $j,k=1,\ldots,N+1.$ For fixed $k,$
  $$\lim\limits_{h\to 0} \sqrt{{\tilde \lambda}_{ik}(h)}=\frac{(2k-1)\pi }{2L}\zeta_i,\quad  i=1,2,\quad k=1,2,\ldots, N.$$
The solutions to \eqref{disc} are given by
\begin{equation*}
\bm \phi(t)= \sum\limits_{j=1}^N  \left[c_{1j}e^{\tilde\lambda_{1j}^+t}{\bm\Psi}^+_{1j}
 + d_{1j}  e^{-\tilde\lambda_{1j}^+t}{\bm\Psi}^-_{1j}+ c_{2j} e^{\tilde\lambda_{2j}^+t}{\bm\Psi}^+_{2j}+d_{2j}e^{-\tilde\lambda_{2j}^+t}{\bm\Psi}^-_{2j}\right].
\end{equation*}
for the initial data
 \begin{eqnarray*}
 \begin{array}{ll}
 \bm \phi_0
 =\sum\limits_{j =1}^N \left( \begin{array}{c}
\frac{1}{\tilde\lambda_{1j}^+}(c_{1j}+ d_{1j})+ \frac{1}{\tilde\lambda_{2j}^+}(c_{2j} + d_{2j}) \\
 \frac{b_1}{\tilde \lambda_{1j}^+} (c_{1j}+ d_{1j})+\frac{b_2}{\tilde\lambda_{2j}^+} (c_{2j} + d_{2j}) \\
(c_{1j}- d_{1j})+(c_{2j} -d_{2j})\\
b_1 (c_{1j}- d_{1j})+b_2(c_{2j} - d_{2j})
 \end{array} \right)\psi_{j}(h)
 \end{array}
\end{eqnarray*}
where $\{c_{ij}, d_{ij}, ~ i=1,2, ~ j\in\mathbb{N}\}$ are complex numbers such that
\begin{align*}
  \| \bm \phi_0\|_{\mathbb{R}^N}^2 \asymp \sum\limits_{j =1}^N \left(|c_{1j}|^2 + |d_{1j}|^2+ |c_{2j}|^2 + |d_{2j}|^2\right), {\rm i.e.}&\\
  \tilde C_1~ \| \bm \phi_0\|_{\mathbb{R}^N}^2 \le \sum\limits_{i=1}^2 \sum\limits_{j =1}^N \left(|c_{ij}|^2 + |d_{ij}|^2\right)\le \tilde C_2 ~\| \bm \phi_0\|_{\mathbb{R}^N}^2
  \end{align*}
with two positive constants $\tilde C_1, \tilde C_2$ which are independent of the particular choice of $\bm \phi_0.$
\end{theorem}

\begin{table}[htb!]
\caption{Construction of the eigenpairs of the matrix $\tilde{\mc A}^{\rm FEM}$ via the eigenpairs of submatrices $\bm C_1^{-1}\bm C_2$ and $\bm M^{-1}\bm A_h$.}
\label{tab:eigenvalues_eigenvectors}
\centering
\begin{tabular}{|c|c|c|c|}
\hline
Matrix & Eigenvalues & Eigenvectors &  \\
\hline
$\bm C_1^{-1}\bm C_2$ & $\zeta_j^2$ & $\begin{bmatrix} 1 \\ b_j \end{bmatrix}$ & $j=1,2$\\
\hline
$\bm M^{-1}\bm A_h$ & $\lambda_k$ & $\bm \psi_k$ & $k=1,2,...,N+1~$ \\
\hline
$\mathbb{A}$ & $\zeta_j^2\lambda_k$ &$\begin{bmatrix} 1 \\ b_j \end{bmatrix} \otimes \bm \psi_k$& $\begin{cases}
k=1,2,...,N+1\\
j=1,2
\end{cases} $
 \\
\hline
$\tilde{\mc A}^{\rm FEM}$ & $\tilde \lambda_{jk}^\mp(h)=\mp i \zeta_j\sqrt{{\lambda}_k }$ & $\begin{bmatrix} \frac{1}{\tilde \lambda_{jk}^\mp(h)} \\ 1 \end{bmatrix} \otimes \begin{bmatrix} 1 \\ b_j \end{bmatrix} \otimes \bm \psi_k$& $\begin{cases}
k=1,2,...,N+1\\
j=1,2
\end{cases} $
\\
\hline
\end{tabular}
\end{table}

\subsection{Lack of Observability as \texorpdfstring{$h\to 0$}{h->0}}
The following lemma provides key identities necessary to prove the main result on the lack of observability.
\begin{lemma} \label{imp1}  For any eigenvector $\begin{bmatrix}
  {\bm \psi}^1&
   {\bm \psi}^2
   \end{bmatrix}^{\rm T}$  of $\mathbb{A}$, the following identities hold:
\begin{equation}\label{10}
\begin{split}
 \alpha_1 \sum\limits_{j=0}^{N+1} \left| \frac{\psi^1_{j+1}-\psi^1_{j}}{h} \right|^2 +\beta \sum\limits_{j=0}^{N+1}  \left| \frac{(\gamma\psi^1_{j+1}-\psi^2_{j+1})-(\gamma\psi^1_{j}-\psi^2_{j})}{h} \right|^2 \\
 \qquad \qquad \qquad={\lambda}_k \left(\sum\limits_{j=0}^{N+1}  \rho |\psi^1_{j}|^2 + \mu |\psi^2_{j}|^2\right),
\end{split}
\end{equation}
and
\begin{equation}\label{13}
\begin{split}
 \frac{{\lambda}_k (2L-h) }{2h} &\left(\rho \left|\psi^1_{N+1}\right|^2 + \mu \left|\psi^2_{N+1}\right|^2\right)\\
 &={\lambda}_k \left[ 2\rho \left(1-\frac{\rho {\lambda}_k  h^2}{12\alpha_1}\right)\sum\limits_{j=0}^{N+1}  |\psi^1_{j}|^2\right.\\
&\quad+  \left. 2\mu \left(1-\frac{\alpha \mu {\lambda}_k  h^2}{12\beta \alpha_1}\right)\sum\limits_{j=0}^N  |\psi^2_{j}|^2-\frac{\rho \mu\gamma{\lambda}_k  h^2}{\alpha_1} \sum\limits_{j=0}^{N+1}  \psi^1_{j} \psi^2_{j}\right].
\end{split}
\end{equation}
\end{lemma}
\begin{proof} The proof follows the same steps as in \cite[Lemma 3.2]{ozer2022uniform}, where the mass matrix $\bm M$ is the identity matrix $\bm I$. Once the eigenvalues of the matrix $\bm A_h$ are replaced by $\lambda_k$, the identities \eqref{10} and \eqref{13} follow directly.
\end{proof}

\begin{theorem}[Lack of Observability]
  \label{thm:observability}
For any $T>0,$
\begin{eqnarray*}
\sup\limits_{ \left( \begin{array}{c}
{\bm v}\\
 {\bm p}
 \end{array} \right)\ ~\rm solves ~(\ref{disc})} \frac{E_h^{\rm FEM}(0)}{\int_0^T \left\{\bm C_1 \begin{bmatrix}
   {\dot {v}}_{N+1} \\
{\dot {p}}_{N+1}     \\
\end{bmatrix}\cdot \begin{bmatrix}
   {\dot {v}}_{N+1} \\
{\dot {p}}_{N+1}     \\
\end{bmatrix}\right\} dt}\to \infty, \text{ as }h\to0.
\end{eqnarray*}
\end{theorem}

\begin{proof} Consider the solution $\left( \begin{array}{c}
\bm v\\
 \bm p
 \end{array} \right)=\left( \begin{array}{c}
\bm \psi_{N+1}\\
b_k {\bm \psi}_{N+1}
 \end{array} \right)$ corresponding to the eigenvalue
 $$
  \tilde \lambda_{N+1}(h)=i \zeta_k \sqrt{{\lambda }_{N+1}}, \quad \text{for } k=1,2.
 $$
By the normalization condition $h\sum\limits_{j=0}^{N+1} |\psi_{N+1,j}|^2 = 1$, equation \eqref{13} gives
  \begin{equation}\label{eq13.3}
    \begin{split}
\frac{\tilde \lambda_{N+1}(2L-h) }{2}\left|\psi_{N+1,N+1}\right|^2 =&E_h^{\rm FEM}(0) \left[ 2\rho \left(\frac{12\alpha_1\zeta_k^2-\rho \lambda_N h^2}{12\alpha_1\zeta_k^2}\right)\right.\\
&+  2\mu b_k^2 \left(1-\frac{\alpha \mu{\lambda}_{N+1} h^2}{12\beta \alpha_1\zeta_k^2}\right) \left.- \frac{\rho \mu\gamma{\lambda}_{N+1} h^2}{\alpha_1\zeta_k^2} b_k \right].
\end{split}
\end{equation}
On the other hand, the observability integral is given by
  \begin{align*}
\int_0^T \left\{C_1 \begin{bmatrix}
   {\dot {v}}_{N+1} \\
{\dot {p}}_{N+1}     \\
\end{bmatrix}\cdot \begin{bmatrix}
   {\dot {v}}_{N+1} \\
{\dot {p}}_{N+1}     \\
\end{bmatrix}\right\} dt&=T\rho\left| \frac{\dot \psi_{N+1,N+1}}{h}\right|+T\mu\left| \frac{\dot \psi_{N+1,N+1}}{h}\right|\\
&=T\tilde \lambda_{N+1}(h) \left(\rho+b_k^2\mu\right)\left| \frac{\psi_{N+1,N+1}}{h}\right|.
\end{align*}
The ratio of the energy to the observability integral becomes
  \begin{equation*}
\begin{split}
&\frac{E_h^{\rm FEM}(0)}{ \int_0^T \left(\rho\left|\dot v_{N+1} \right|^2+\mu\left|\dot p_{N+1} \right|^2\right) dt}  =(2L-h)/ \\
&\hspace*{3cm} \left\{2T(\rho+\mu b_k^2)\left[ 2\rho \left(\frac{12\alpha_1\zeta_k^2-\rho {\lambda}_{N+1}h^2}{12\alpha_1\zeta_k^2}\right)\right.\right.\\
&\hspace*{3cm} \qquad \left.\left. +  2\mu b_k^2 \left(1-\frac{\alpha \mu {\lambda}_{N+1} h^2}{12\beta \alpha_1\zeta_k^2}\right)- \frac{\rho \mu\gamma {\lambda}_{N+1} h^2}{\alpha_1\zeta_k^2} b_k \right]\right\}.
\end{split}
\end{equation*}
Since $\lambda_{N+1} \to 12$ as $h \to 0$, equation \eqref{eq13.3} reduces to
   \begin{equation}\label{eq13.4}
\begin{split}
\frac{E_h^{\rm FEM}(0)}{ \int_0^T \left(\rho\left|\dot v_{N+1} \right|^2+\mu\left|\dot p_{N+1} \right|^2\right) dt} &=L\left/ \left\{T(\rho+\mu b_k^2)\left[ 2\rho \left(\frac{\alpha_1\zeta_k^2-\rho }{\alpha_1\zeta_k^2}\right)\right.\right.\right.\\
&\qquad\left.\left. +  2\mu b_k^2 \left(1-\frac{\alpha \mu }{\beta \alpha_1\zeta_k^2}\right)- \frac{12\rho \mu\gamma}{\alpha_1\zeta_k^2} b_k \right]\right\}.
\end{split}
\end{equation}
By utilizing \eqref{b-quad}, factoring out $\frac{2\mu}{\beta(b_k \gamma \mu + \rho)}$ in the denominator of \eqref{eq13.4}, and analyzing the limit as $h \to 0$, the ratio diverges
\begin{align*}
&\frac{E_h^{\rm FEM}(0)}{ \int_0^T \left(\rho\left|\dot v_N \right|^2+\mu\left|\dot p_N \right|^2\right) dt}  \to \frac{L}{ T(\rho+\mu b_k^2)\left[  2\mu b_k^2 \left(1-\frac{\alpha \mu \zeta_k^2 }{\beta \alpha_1}\right)- \frac{2\rho \mu\gamma \zeta_k^2 }{\alpha_1} b_k \right]}\\
&\hspace*{1cm} \to \frac{\beta (b_k \gamma \mu +\rho)}{2\mu b_k} \frac{L}{T(\rho+\mu b_k^2)\left[  b_k \left(\beta (b_k \gamma\mu +\rho)-\alpha \mu \right)- \rho\gamma  \beta \right]}\\
& \hspace*{1cm}\to  \frac{\beta (b_k \gamma \mu +\rho)}{2\mu b_k}\frac{L}{{T(\rho+\mu b_k^2)\left[  \beta \gamma\mu b_k^2 + \rho \gamma b_k - \alpha \mu b_k - \rho \gamma \beta \right]}}\\
& \hspace*{1cm}\to \frac{\beta (b_k \gamma \mu +\rho)}{2\mu b_k}\frac{L}{{T(\rho+\mu b_k^2)\left[ \gamma \beta \mu b_k^2  - (\alpha \mu - \rho \beta)b_k -\rho \gamma \beta \right]}} \to \infty,
\end{align*}
where   $\beta \gamma\mu b_k^2 -(\alpha \mu - \beta \rho)b_k - \rho\gamma\beta=0$ by \eqref{b-quad}.
\end{proof}

This result shows that the discrete energy cannot be recovered from observations at the right boundary $x = L$. As a result, the solutions of the closed-loop system \eqref{FEM-pseudo} cannot be exponentially stable uniformly as $h \to 0$. The outcome aligns with \cite{ozer2023revisiting}, where observability inequalities and exponential stability fail uniformly as $h \to 0$ in standard Finite Difference approximations of \eqref{disc}.

\subsection{Direct Fourier Filtering for Preserving Exact Observability as  \texorpdfstring{$h\to 0$}{h->0}}
\label{sec:filtering}

The exact observability result for the PDE model was studied in \cite{ozer2022uniform}. However, when discretizing the PDE model using Finite Element methods, the lack of observability demonstrated in \cref{thm:observability} becomes a significant limitation as $h \to 0$. This issue arises primarily from the high-frequency eigenvalues of the matrix $\bm A_h$, where ${\lambda}_{N+1} \to 12$ as $h \to 0$. A well-established remedy for this issue is the direct Fourier filtering technique, which has been successfully implemented in \cite{ozer2022uniform, ozer2023revisiting} for Finite Difference approximations. These studies demonstrated that the observability inequality and exponential stability results could be retained uniformly as $h \to 0$ using this technique.

The goal here is to extend the direct Fourier filtering approach to the Finite Element approximations to recover the exact observability of the PDE model uniformly as $h \to 0$.

To define the direct Fourier filtering, let $j^*$ be the filtering parameter, as introduced in \cite{ozer2022uniform}. The class of filtered solutions for the Finite Element approximation \eqref{FEM} is represented by
\begin{eqnarray}
  \label{eq:filtered}
	\begin{array}{ll} {\mathcal C}_h(j^*):=\sum\limits_{j=1}^{N+1-j^*}  \left[c_{1j}e^{\lambda^+_{1j}t}{\bm\Psi}^{+}_{1j}
		+ d_{1j} e^{{{\lambda}}^-_{1j}t} {\bm{\Psi}^{-}_{1j}} + c_{2j} e^{\lambda^+_{2j}t}{\bm\Psi}^{+}_{2j}+d_{2j}e^{-{\lambda}^-_{2j}t}{\bm{ \Psi}^{-}_{2j}}\right] \end{array}
\end{eqnarray}
where $j^*$ specifies the number of eigenvalues filtered from each signed half of each branch. From the total $4N+4$ eigenpairs of the matrix $\tilde{\mc A}^{\rm FEM}$, this filtering method removes $4j^*$ eigenpairs from the solutions.

The direct Fourier filtering technique targets the high-frequency eigenvalues responsible for the lack of observability. As shown in \cref{tab:eigenvalues_eigenvectors}, regardless of the positive/negative sign or the branch selection ($\zeta_1$ or $\zeta_2$), the eigenvalues $\tilde \lambda_{j(N+1)}^\mp$ always include the same ${\lambda}_{N+1}$ term. Consequently, the filtering is applied independently to each signed half of each branch.

A critical step in this process is the proper separation of eigenpairs associated with each branch, which ensures the accurate application of the filtering. An algorithm addressing this separation is presented in \cref{rmk:filtering}, and a detailed numerical application is provided in \cref{Simu}. Although the Fourier filtering approach remedies the lack of observability, a direct proof of exponential stability via a Lyapunov-based method is beyond the scope of this paper. Such a proof would require novel spectral estimates and deeper analytical insights into the system's dynamics.

{
\section{Proposed Model Reduction by Order Reduction-based Finite Differences (ORFD)}
\label{ModelRed}
As demonstrated in the previous section, using standard Finite Element or Finite Difference \cite{ozer2022uniform} approximations for \eqref{eq1} leads to difficulties in preserving exponential stability uniformly as $h \to 0$. To address this issue, we adopt a novel order-reduction methodology inspired by the work in \cite{aydinhaider2024, aydin2023novel, liu2020new} to construct a new model reduction.

We begin by introducing the following variables for reducing the order of the system \eqref{eq1}
\begin{equation}
\label{sub}
u^1(x,t)=v_x(x,t),~ u^2(x,t)=p_x(x,t),~
w^1(x,t)=\dot v(x,t), ~ w^2(x,t)=\dot p(x,t).
\end{equation}
Substituting \eqref{sub} into \eqref{eq1} and appending the resulting equations, the PDE system \eqref{eq1} is reformulated as
\begin{equation}
\label{or1}
\begin{cases}
\rho {\dot w}^1-\alpha u^1_{x}+ \gamma \beta u^2_{x}=0,\\
\mu {\dot w}^2-\beta u^2_{x}+\gamma \beta u^1_{x}=0,&\\
{\dot u}^1-w^1_x=0,\\
{\dot u}^2-w^2_x=0, &(x,t) \in (0,L) \times  \mathbb{R}^+ \\
w^1(0,t)=w^2(0,t)=0,&\\
\alpha u^1(L,t)-\gamma \beta u^2(L,t)= -\xi_1 w^1(L,t),&\\
\beta u^2(L)-\gamma\beta u^1(L)=-\xi_2 w^2(L,t),&t\in \mathbb{R}^+\\
(u^1,w^1,u^2,w^2)(x,0) = (\frac{d}{dx}[v_0], v_1,\frac{d}{dx}[p_0], p_1)(x),& x\in[0,L].
\end{cases}
\end{equation}
}
Consider the same uniform discretization of $[0,L]$.  We introduce grid functions, $U_j^1(t),U_j^2(t),W_j^1(t),W_j^2(t)$ for all  $j\in \{0,1,2,3,...,N+1\}$ as follows
\begin{equation*}
U_j^1(t)\approx u^1(x_j,t),~U_j^2(t) \approx u^2(x_j,t),~
W_j^1(t) \approx w^1(x_j,t),~W_j^2(t) \approx w^2(x_j,t).
\end{equation*}
The averaging and difference operators are defined at the in-between nodes $x_{j+\frac{1}{2}}=\frac{x_j+x_{j+1}}{2}=\left(j+\frac{1}{2}\right)h$ for $j=0,1,\ldots, N.$ The following approximations at these nodes are defined
\begin{eqnarray*}
U_{j+\frac{1}{2}}:=\frac{U_j+U_{j+1}}{2}, \qquad \delta_x[U_{j+\frac{1}{2}}]:=\frac{U_{j+1}-U_j}{h} .
\end{eqnarray*}

Next, we approximate the equations $\eqref{or1}_1$-$\eqref{or1}_4$ at each midpoint $\left\{(x_{j+\frac{1}{2}},t)\right\}_{j=1}^{N}$ to obtain
\begin{equation}
	\label{ORFD}
\begin{cases}
	{\dot u}^{1}(x_{j+\frac{1}{2}}) - \mathcal {D}_x[w^1](x_{j+\frac{1}{2}}) =0,\\
{\dot u}^{2}(x_{j+\frac{1}{2}})-\mathcal {D}_x[w^2](x_{j+\frac{1}{2}}) =0,&\\
\rho {\dot w}^{1}(x_{j+\frac{1}{2}})-\alpha  \mathcal {D}_x[u^1](x_{j+\frac{1}{2}}) + \gamma \beta\ \mathcal {D}_x [u^2](x_{j+\frac{1}{2}}) = 0,&\\
\mu {\dot w}^{2}(x_{j+\frac{1}{2}}) -\beta \mathcal {D}_x [u^2](x_{j+\frac{1}{2}}) +\gamma \beta \mathcal {D}_x[u^1](x_{j+\frac{1}{2}}) =0,\qquad j=1,\ldots, N,&\\
	\end{cases}
\end{equation}
where $\mathcal {D}_x=\frac{\partial }{\partial x},$ and ${\mathcal D}_x^0=\mathcal I$ is the identity operator. Using Taylor's theorem, if $f (x) \in C^2[0,L]$ and $g(x) \in C^3[0,L]$, the following approximations hold
\begin{equation}\label{ornek1}
\begin{split}
 f_{j+\frac{1}{2}}=&{\mathcal D}_x^0[f]\left(x_{j+\frac{1}{2}}\right)+\frac{h^2}{8}\mathcal {D}_{xx}[{f}]\left(\tilde \xi_1 \right), \\
\delta _x\left[g_{j+\frac{1}{2}}\right]=&\mathcal {D}_x[g]\left(x_{j+\frac{1}{2}}\right)+\frac{h^2}{24} \mathcal {D}_{xxx}[g](\tilde \xi_2)
\end{split}
\end{equation}
 for some $\tilde \xi_1,\tilde \xi_2\in \left[x_j,x_{j+1}\right]$.
Replacing the operator $ {\mathcal D}x$ by $\delta_x$ in \eqref{ORFD}, approximating the initial and boundary conditions  $\eqref{or1}_5$-$\eqref{or1}_8$ at the corresponding nodes, and gathering the higher-order terms on the right-hand side, we obtain the following equivalent system
\begin{equation}
\label{discmodelwremainder}
\begin{cases}
{\dot U}^{1}_{j+\frac{1}{2}} - \delta_x[W^1_{j+\frac{1}{2}}] =r_{j+\frac{1}{2}}^1(t),&\\
{\dot U}^{2}_{j+\frac{1}{2}}- \delta_x[W^2_{j+\frac{1}{2}}] =r_{j+\frac{1}{2}}^2(t),&\\
\rho {\dot W}^{1}_{j+\frac{1}{2}}-\alpha\, \delta_x[U^1_{j+\frac{1}{2}}] + \gamma \beta\; \delta_x [U^2_{j+\frac{1}{2}}] = s_{j+\frac{1}{2}}^1(t),&\\
\mu {\dot W}^{2'}_{j+\frac{1}{2}} -\beta\,\delta_x [U^2_{j+\frac{1}{2}}] +\gamma \beta\,  \delta_x[U^1_{j+\frac{1}{2}}] =s_{j+\frac{1}{2}}^2(t),& j=1,\ldots, N,\\
W^1_0=W^2_0=0, \\
\alpha U^1_{N+1}-\gamma \beta U^2_{N+1}= -k_1 W^1_{N+1},&\\
\beta U^2_{N+1}-\gamma\beta U^1_{N+1}=-k_2 W^2_{N+1}, &  t\in \mathbb{R}^+,\\
(U^1_j, U^2_j, W^1_j, W^2_j)(0) (  {\mathcal D}_{x}[w_0],   {\mathcal D}_{x}[p_0], w_1,p_1)(x_j),& j=0,\ldots, N+1,
\end{cases}
\end{equation}
for $j=0,\ldots, N+1,$ where  the ``residual" terms $r^1,s^1,r^2,s^2$ are bounded by \eqref{ornek1} by the assumption that $ {u}^k (x), w^k\in C^3[0,L]$ for $k=1,2.$

The residual terms in \eqref{discmodelwremainder} are explicitly given by
\begin{equation}
\label{remainders}
\begin{cases}
s_{j+\frac{1}{2}}^1(t) = h^2 \left( \frac{\rho }{8} {\mathcal D}_{xx}[w^1](\xi_1) - \frac{\alpha}{24} {\mathcal D}_{xxx}[u^1](\xi_2)+\frac{ \gamma \beta}{24}  {\mathcal D}_{xxx}[u^2](\xi_3) \right), \\
s_{j+\frac{1}{2}}^2(t)=  h^2 \left( \frac{\mu }{8}  {\mathcal D}_{xx}[w^2](\xi_4) - \frac{\beta}{24}  {\mathcal D}_{xxx}[u^2](\xi_3)+\frac{ \gamma \beta}{24}  {\mathcal D}_{xxx}[u^1](\xi_2) \right),\\
 r_{j+\frac{1}{2}}^1(t)= h^2\left(\frac{1}{8}  {\mathcal D}_{xx}[u^1](\xi_5)-\frac{1}{24}  {\mathcal D}_{xxx}[w^1](\xi_6)\right),  \\
r_{j+\frac{1}{2}}^2(t)= h^2\left(\frac{1}{8}  {\mathcal D}_{xx}[u^2](\xi_7)-\frac{1}{24}  {\mathcal D}_{xxx}[w^2](\xi_8)\right),
\end{cases}
\end{equation}
 for some  $\{\xi_1,\xi_2,...,\xi_8\} \in [x_j,x_{j+1}].$

We assume that the initial data $( {\mathcal D}_{x}[w_0],  {\mathcal D}_{x}[p_0], w_1,p_1)(x)$ is smooth enough to allow for evaluations at each grid point $x_j$, for $j = 0,1,\dots,N+1$.

By neglecting the infinitesimal terms $r^1, s^1, r^2, s^2$ in \eqref{discmodelwremainder} and replacing $(U^1_j, U^2_j, W^1_j, W^2_j)$ with $(u^1_j, u^2_j, w^1_j, w^2_j),$ we obtain the semi-discretized finite difference approximation for \eqref{or1}
\begin{equation}
	\label{ORFD2}
\begin{cases}
{\dot u}^{1}_{j+\frac{1}{2}} - \delta_x[w^1_{j+\frac{1}{2}}] =0,\\
{\dot u}^{2}_{j+\frac{1}{2}}- \delta_x[w^2_{j+\frac{1}{2}}] =0,&\\
\rho {\dot w}^{1}_{j+\frac{1}{2}}-\alpha\, \delta_x[u^1_{j+\frac{1}{2}}] + \gamma \beta\; \delta_x [u^2_{j+\frac{1}{2}}] = 0,&\\
\mu {\dot w}^{2'}_{j+\frac{1}{2}} -\beta\,\delta_x [u^2_{j+\frac{1}{2}}] +\gamma \beta\,  \delta_x[u^1_{j+\frac{1}{2}}] =0,& j=1,\ldots, N,\\
u^1_0=u^2_0=0,\\
 \alpha u^1_{N+1}-\gamma \beta u^2_{N+1}= -k_1 w^1_{N+1},&\\
\beta u^2_{N+1}-\gamma\beta u^1_{N+1}=-k_2 w^2_{N+1}, & t\in \mathbb{R}^+,\\
(u^1, u^2, w^1, w^2)_j(0)=  (u^1, u^2, w^1, w^2)_j^0, & j=0,\ldots, N+1,
\end{cases}
\end{equation}
\begin{lemma}\cite[Lemma 2.2]{liu2020new}\label{summa}
For any grid functions $\bm{U}, \bm{V}, \bm{W}$ defined on the mesh points $\{x_j\}_{j=0}^{N+1}$, the following summation by parts formulas hold:
\begin{equation*}
\begin{array}{ll}
    (i) ~~ h \sum\limits _{j=0}^N \delta_x U_{j+\frac{1}{2}} V_{j+\frac{1}{2}} + h \sum\limits _{j=0}^N \delta_x V_{j+\frac{1}{2}} U_{j+\frac{1}{2}} + U_0V_0 = U_{N+1}V_{N+1}, \\[10pt]

    (ii)~~ h \sum\limits _{j=0}^N \delta_x (U_{j+\frac{1}{2}} V_{j+\frac{1}{2}} W_{j+\frac{1}{2}}) + U_0V_0W_0
    + h \sum\limits_{j=0}^N U_{j+\frac{1}{2}}  \delta_x V_{j+\frac{1}{2}} W_{j+\frac{1}{2}} \\[5pt]
    \quad + h \sum\limits_{j=0}^N U_{j+\frac{1}{2}} V_{j+\frac{1}{2}}  \delta_x W_{j+\frac{1}{2}} = U_{N+1}V_{N+1}W_{N+1}  \\
\quad    - \frac{h}{4} \sum\limits _{j=0}^N (U_{j+1}-U_j) (V_{j+1}-V_j) (W_{j+1}-W_j).
\end{array}
\end{equation*}
\end{lemma}

The discretized energy for \eqref{ORFD2} can be redefined as the following
\begin{equation}
\label{enpde2}
\begin{array}{ll}
	 E^{\rm ORFD}_h(t):=\frac{h}{2} \sum\limits_{j=0}^N
\left(\rho     |w^1_{j+\frac{1}{2}}(t)|^2  +\mu    |w^2_{j+\frac{1}{2}}|^2 \right.\\
\left.\qquad \qquad \qquad  +\beta  |\gamma  u^1_{j+\frac{1}{2}}-u^2_{j+\frac{1}{2}}| ^2
+\alpha_1| u^1_{j+\frac{1}{2}}| ^2\right).
\end{array}
\end{equation}

\begin{lemma} \label{lemm1}The energy $E^{\rm ORFD}_h(t)$ is dissipative, i.e.
\begin{equation}\label{lemres-b}
    \dot E^{\rm ORFD}_h(t)= -k_1\left|w^1 _{N+1}\right|^2- k_2\left|w^2 _{N+1}\right|\le 0.
\end{equation}
\end{lemma}
\begin{proof} By taking the derivative of the energy $E^{\rm ORFD}_h(t)$ along the solutions of \eqref{ORFD2}, we obtain
\begin{align*}
   \dot E^{\rm ORFD}_h(t)
   =\hspace*{-0.02cm}&\sum\limits _{j=0}^N  w^1_{j+\frac{1}{2}} \left(\alpha\ \delta_x[u^1_{j+\frac{1}{2}}] - \gamma \beta \delta_x [u^2_{j+\frac{1}{2}}] \right)+ w^2_{j+\frac{1}{2}} \beta\left( \delta_x [u^2_{j+\frac{1}{2}}] -\gamma   \delta_x[u^1_{j+\frac{1}{2}}]\right)\\
  &+ \beta \left(\gamma u^1_{j+\frac{1}{2}}-u^2_{j+\frac{1}{2}}\right) \left(\gamma \delta_x[w^1_{j+\frac{1}{2}}]-\delta_x[w^2_{j+\frac{1}{2}}]\right) + \alpha_1 \rho u^1_{j+\frac{1}{2}} \delta_x[w^1_{j+\frac{1}{2}}].
\end{align*}
Now, using the first result in Lemma \ref{summa}, we obtain
\begin{align*}
 \dot E^{\rm ORFD}_h(t)=&
   \sum\limits _{j=0}^N  -\alpha_1 \delta_x[w^1_{j+}]  u^1_{j+\frac{1}{2}} - \beta \left(\gamma u^1_{j+\frac{1}{2}}-u^2_{j+\frac{1}{2}}\right) \left(\gamma \delta_x[w^1_{j+\frac{1}{2}}]-\delta_x[w^2_{j+\frac{1}{2}}]\right)\\
&+ \beta \left(\gamma u^1_{j+\frac{1}{2}}-u^2_{j+\frac{1}{2}}\right) \left(\gamma \delta_x[w^1_{j+\frac{1}{2}}]-\delta_x[w^2_{j+\frac{1}{2}}]\right)\\
   &+ \alpha_1 \rho u^1_{j+\frac{1}{2}} \delta_x[w^1_{j+\frac{1}{2}}]  -w^2_{0} \left(\beta u^2_{0} -\gamma \beta  u^1_{0}\right) \\
& \hspace*{-1.5cm}+ w^2_{N+1} \left(\beta u^2_{N+1} -\gamma \beta  u^1_{N+1}\right)    -w^1_0 \left(\alpha u^1_{0} - \gamma \beta u^2_{0} \right)+w^1_{N+1} \left(\alpha u^1_{N+1} - \gamma \beta u^2_{N+1} \right).
\end{align*}
Finally, applying the boundary conditions completes the proof of \eqref{lemres-b}.
\end{proof}

\begin{remark}
  The system \eqref{ORFD2} is formulated in terms of lower-order variables. However, it can be expressed in terms of the original variables $v$ and $p$ using the following approximations
  \begin{align*}
  u^1_{j+\frac{1}{2}}\approx \delta_x[v_{j+\frac{1}{2}}], \quad u^2_{j+\frac{1}{2}}\approx \delta_x[p_{j+\frac{1}{2}}],  \quad w^1_{j+\frac{1}{2}}\approx \dot v_{j+\frac{1}{2}}, \quad w^2_{j+\frac{1}{2}}\approx \dot p_{j+\frac{1}{2}}.
  \end{align*}
Consequently, the system \eqref{ORFD2} can be rewritten in terms of the original variables in a matrix form analogous to the FEM approximation in \cref{FEMSplines}
\begin{equation}
    \label{eq:ORFD_matrix}
\begin{cases}
  (\bm C_1 \otimes \bm M_{\rm ORFD}) \begin{bmatrix}
       {\ddot{\bm{v}}}   \\
     {\ddot{\bm{p}}}        \\
  \end{bmatrix} +\left(\bm C_2  \otimes \bm A_h\right) \begin{bmatrix}
       {{\bm{v}}}   \\
     {{\bm{p}}}        \\
  \end{bmatrix}	 + (\bm C_3 \otimes \bm B) \begin{bmatrix}
       {\dot {\bm{v}}}   \\
     {\dot {\bm{p}}}        \\
  \end{bmatrix} =\bm 0,  & t\in \mathbb{R}^+,\\
  (\bm v, \dot{\bm{v}},\bm p,\dot{\bm p})_i(0)=(v^0,v^1,p^0, p^1)(x_i),  & \hspace*{-1cm}i=0,...,N+1,
    \end{cases}
  \end{equation}
where $\bm M_{\rm OR}$ is the $(N+1) \times (N+1)$ mass matrix, given by
\begin{equation} \label{eq:massOR}
  \bm M_{\rm OR}:=\frac{1}{4}\begin{bmatrix}
		2&1&0&\dots&\dots&\dots&0\\
		1&2&1&0&\dots&\dots&0\\
		&\ddots&\ddots&\ddots&\ddots&\ddots&\\
		0&\dots&\dots&0&1&2&1\\
		0&\dots&\dots&\dots&0&1&1\\
	\end{bmatrix}.
\end{equation}
Furthermore, the approximated energy for the system \eqref{eq:ORFD_matrix} is given by
\begin{equation}\label{ORFFDEnergy}
    E^{\rm ORFD}_h(t):=\frac{h}{2}\left\langle (\bm C_1 \otimes \bm M_{\rm OR})  \begin{bmatrix}
			{\dot {\vec v}} \\
			{\dot {\vec p}}
		\end{bmatrix},\begin{bmatrix}
			{\dot {\vec v}} \\
			{\dot {\vec p}}
		\end{bmatrix}\right\rangle
 +\frac{h}{2}\left\langle (\bm C_2 \otimes \bm A_h)  \begin{bmatrix}
			{ {\vec v}} \\
			{ {\vec p}}
		\end{bmatrix},\begin{bmatrix}
			{ {\vec v}} \\
			{ {\vec p}}
		\end{bmatrix}\right\rangle.
\end{equation}
\end{remark}

  \subsection{Exponential Stability by the Lyapunov Approach}
  To establish the exponential stability of the system, we introduce a Lyapunov functional that combines the discrete energy with an auxiliary term designed to capture the system's decay properties.

For $\delta>0$, define the Lyapunov functional as
  \begin{equation*}
  L_h(t):=E^{\rm ORFD}_h(t)+ \delta \phi_h(t)
\end{equation*}
where the auxiliary functional $\phi_h(t)$ is defined by
  \begin{equation}
   \phi_h(t):=h\sum\limits _{j=0}^N x_{j+\frac{1}{2}} (\rho u^1_{j+\frac{1}{2}}w^1_{j+\frac{1}{2}}+\mu u^2_{j+\frac{1}{2}} w^2_{j+\frac{1}{2}}).
\end{equation}
\begin{lemma}  \label{lemm2}
For $\delta<\frac{1}{L\eta }, $ the Lyapunov function $L_h(t)$ and the energy $E^{\rm ORFD}_h(t)$ in \eqref{enpde2} are equivalent, i.e. there exists $C_1,C_2>0$ such that
\begin{equation}\label{lemres-c}
    C_1 E^{\rm ORFD}_h(t) \le L_h(t)\le C_2 E^{\rm ORFD}_h(t)
\end{equation}
where $C_1=1-\delta L \eta, $ $C_2=1+\delta L \eta,$ and $\eta$ is defined by \eqref{eta}.
\end{lemma}
\begin{proof}
Using Hölder's, Cauchy-Schwarz's, and Minkowski's inequalities, we obtain
\begin{align*}
&|\phi_h(t)| \le
 h L \sqrt{ \frac{\rho}{\alpha_1}}\sum\limits _{j=0}^N \sqrt{\rho} u^1_{j+\frac{1}{2}} \sqrt{\alpha_1}w^1_{j+\frac{1}{2}}+h L \sqrt{ \frac{\mu}{\beta}}\sum\limits _{j=0}^N \sqrt{\mu} u^2_{j+\frac{1}{2}} \sqrt{\beta}w^2_{j+\frac{1}{2}}\\
&\qquad \quad \le h L  \sqrt{ \frac{\rho}{\alpha_1}} \left[\left(\sum\limits _{j=0}^N \rho |u^1_{j+\frac{1}{2}} |^2\right) \left(\sum\limits _{j=0}^N \alpha_1 |w^1_{j+\frac{1}{2}} |^2\right)\right]^{\frac{1}{2}}\\
&+ h L  \left(\sum\limits _{j=0}^N  \mu |u^2_{j+\frac{1}{2}}|^2\right)^{\frac{1}{2}}\left\{ \mu\beta\left| \sum\limits _{j=0}^N w^2_{j+\frac{1}{2}}-\gamma w^1_{j+\frac{1}{2}}\right|^2 + \left(\frac{\mu\gamma^2}{ \alpha_1}  \sum\limits _{j=0}^N \alpha_1w^1_{j+\frac{1}{2}}\right)^{\frac{1}{2}}\right\}\\
&\qquad \quad =\frac{h L}{2}\left\{  \sqrt{ \frac{\rho}{\alpha_1}}\sum\limits _{j=0}^N \rho |u^1_{j+\frac{1}{2}} |^2 \right. +\sqrt{\frac{\mu}{\beta}}\sum\limits _{j=0}^N \beta   |w^2_{j+\frac{1}{2}}-\gamma w^1_{j+\frac{1}{2}}|^2\\
& \qquad \qquad  +\left( \sqrt{ \frac{\rho}{\alpha_1}}+\sqrt{\frac{\mu\gamma^2}{ \alpha_1} } \right) \sum\limits _{j=0}^N \alpha_1 |w^1_{j+\frac{1}{2}} |^2\left.+\left( \sqrt{ \frac{\mu}{\beta}}+\sqrt{\frac{\mu\gamma^2}{ \alpha_1} } \right) \sum\limits _{j=0}^N \mu |u^2_{j+\frac{1}{2}} |^2\right\}.
\end{align*}
Thus, \eqref{lemres-c} follows from the definition of $\eta$ in \eqref{eta}.
\end{proof}

\begin{lemma}\label{lemm3}
The auxiliary function $\phi_h(t)$ satisfies the following inequality
\begin{equation}\label{lemres-d}
\dot \phi_h(t)\hspace*{-0.1cm}\le \hspace*{-0.1cm}-E^{\rm ORFD}_h(t)+  \frac{L}{2} \left\{\left(\rho + \frac{k_1^2}{\alpha_1}\right) |w^1_{N+1}|^2  + \left(\mu+\frac{(\alpha+\gamma^2 \beta)k_2^2}{2\alpha_1\beta}\right) |w^2_{N+1}|^2\right\}.
\end{equation}
\end{lemma}
\begin{proof}
A direct calculation gives
\begin{align*}
\dot \phi_h(t)=&h\sum\limits _{j=0}^N  x_{j+\frac{1}{2}}\left\{ \rho ({\dot u}^1_{j+\frac{1}{2}}w^1_{j+\frac{1}{2}}+u^1_{j+\frac{1}{2}}{\dot w}^1_{j+\frac{1}{2}} )  +\mu ({\dot u}^2_{j+\frac{1}{2}} w^2_{j+\frac{1}{2}}+u^2_{j+\frac{1}{2}} {\dot w}^2_{j+\frac{1}{2}}) \right\}\\
=&h\sum\limits _{j=0}^N x_{j+\frac{1}{2}} \left\{\rho w^1_{j+\frac{1}{2}}  \delta_x[w^1_{j+\frac{1}{2}}] +\mu w^2_{j+\frac{1}{2}}  \delta_x[w^2_{j+\frac{1}{2}}]  \right.\\
&+u^1_{j+\frac{1}{2}} \left(\alpha \delta_x[u^1_{j+\frac{1}{2}}] - \gamma \beta \delta_x [u^2_{j+\frac{1}{2}}] \right) \left. + u^2_{j+\frac{1}{2}} \left(\beta \delta_x[u^2_{j+\frac{1}{2}}] - \gamma \beta \delta_x [u^1_{j+\frac{1}{2}}] \right) \right\}.
\end{align*}
Now, applying the second summation by parts formula from Lemma \ref{summa}, we obtain
\begin{align*}
\dot \phi_h(t)=&L\left(\rho |w^1_{N+1}|^2+\mu |w^2_{N+1}|^2\right) \\
&+ h\sum\limits _{j=0}^N \left\{ -\rho|w^1_{j+\frac{1}{2}}|^2-\mu |w^2_{j+\frac{1}{2}}|^2  -\frac{ \rho}{4}  |w^1_{j+1}-w^1_{j}|^2-\frac{ \mu}{4}   |w^2_{j+1}-w^2_{j}|^2 \right\}\\
& + \sum\limits _{j=0}^N x_{j+\frac{1}{2}} \left\{ - \rho \delta_x[w^1_{j+\frac{1}{2}}] w^1_{j+\frac{1}{2}} \right.-\mu    \delta_x[w^2_{j+\frac{1}{2}}] w^2_{j+\frac{1}{2}}+ \alpha_1 u^1_{j+\frac{1}{2}}  \delta_x[u^1_{j+\frac{1}{2}}] \\
& \left.+ \beta (\gamma u^1_{j+\frac{1}{2}}-u^2_{j+\frac{1}{2}}) \delta_x[\gamma u^1_{j+\frac{1}{2}}-u^2_{j+\frac{1}{2}}] \right\},
\end{align*}
and therefore,
\begin{align*}
\dot \phi_h(t)=& L\left\{ \rho |w^1_{N+1}|^2+\mu |w^2_{N+1}|^2+ \alpha_1 |u^1_{N+1}|^2+ \beta |\gamma u^1_{N+1}-u^2_{N+1}|^2\right\}\\
& + \sum\limits _{j=0}^N \left\{ -\rho |w^1_{j+\frac{1}{2}}|^2-\mu |w^2_{j+\frac{1}{2}}|^2 - \alpha_1  |u^1_{j+\frac{1}{2}}|^2-h\beta  |\gamma u^1_{j+\frac{1}{2}}-\gamma u^2_{j+\frac{1}{2}}|^2 \right\}\\
& + h  \sum\limits _{j=0}^N x_{j+\frac{1}{2}} \left\{-\rho  \delta_x[w^1_{j+\frac{1}{2}}] w^1_{j+\frac{1}{2}} \right.- \mu \delta_x[w^2_{j+\frac{1}{2}}] w^2_{j+\frac{1}{2}} -  \alpha_1 \delta_x[u^1_{j+\frac{1}{2}}]  u^1_{j+\frac{1}{2}} \\
&  \left.-  \beta \delta_x \left[\gamma u^1_{j+\frac{1}{2}}-u^1_{j+\frac{1}{2}}\right] \left(\gamma u^1_{j+\frac{1}{2}}-u^1_{j+\frac{1}{2}}\right)\right\}\\
& + \frac{h }{4}  \sum\limits _{j=0}^N\left\{ -\rho|w^1_{j+1}-w^1_{j}|^2-\mu |w^2_{j+1}-w^2_{j}|^2 \right.\\
 & \left. -\alpha_1|u^1_{j+1}-u^1_{j}|^2-\beta |(\gamma u^1-u^2)_{j+1}-(\gamma u^1-u^2)_{j}|^2 \right\}.
\end{align*}
Using the definition of $\phi_h(t)$
\begin{align*}
\dot \phi_h(t)=& - \frac{d\phi_h(t)}{dt}-2 E^{\rm ORFD}_h(t)+L\left\{ \rho |w^1_{N+1}|^2 + \mu |w^2_{N+1}|^2+ \alpha_1 |u^1_{N+1}|^2\right.  \\
& \left.+ \beta |\gamma u^1_{N+1}-u^2_{N+1}|^2\right\} + \frac{h}{4}  \sum\limits _{j=0}^N \left\{ -\rho |w^1_{j+1}-w^1_{j}|^2- \mu |w^2_{j+1}-w^2_{j}|^2 \right. \\
 &\left. -\alpha_1 |u^1_{j+1}-u^1_{j}|^2-\beta |(\gamma u^1-u^2)_{j+1}-(\gamma u^1-u^2)_{j}|^2\right\},
\end{align*}
and discarding the non-positive terms and using the boundary conditions \eqref{ORFD2}
\begin{align*}
\dot \phi_h(t)\le& - E^{\rm ORFD}_h(t)+\frac{L}{2} \left\{ \rho |w^1_{N+1}|^2 + \mu|w^2_{N+1}|^2 \right. \\
&\left. +\alpha_1 \left|-\frac{k_1}{\alpha_1}w^1_{N+1}-\frac{\gamma k_2}{\alpha_1} w^2_{N+1}\right|^2+ \beta\left|\frac{k_2}{\beta} w^2_{N+1}\right|^2 \right\}.
\end{align*}
Finally, applying the Cauchy-Schwarz inequality completes the proof.
\end{proof}

\begin{theorem}\label{mainexp} [Exponential Stability] For any $k_1,k_2>0,$ and
\begin{equation}\label{deltacap}
    \delta<\frac{1}{2L}\min \left(\frac{1}{\eta},\frac{2k_1\alpha_1}{\alpha_1\rho + k_1^2}, \frac{4k_2\beta\alpha_1}{2\alpha_1\beta\mu+(\alpha+\gamma^2 \beta)k_2^2}\right),
\end{equation}
the energy $E^{\rm ORFD}_h(t)$ in \eqref{enpde2} decays exponentially, i.e.
\begin{equation}\label{mainthm1}
\dot E^{\rm ORFD}_h(t)\le \frac{1+\delta L\eta}{1-\delta L\eta} e^{-\delta (1-\delta L\eta) t} E^{\rm ORFD}_h(0).
\end{equation}
\end{theorem}
\begin{proof}
By Lemmas \ref{lemm1} and  \ref{lemm3}, we obtain
\begin{align*}
  \dot L_h(t)=&  \dot E^{\rm ORFD}_h(t)  + \delta  \dot \phi_h(t) \\
   =&-\delta E^{\rm ORFD}_h(t)-\left(k_1-\frac{\delta L}{2}\left(\rho + \frac{k_1^2}{\alpha_1}\right) |w^1_{N+1}|^2\right)\\
  &  -\left(k_2-\frac{\delta L}{2} \left(\mu+\frac{(\alpha+\gamma^2 \beta)k_2^2}{2\alpha_1\beta}\right) |w^1_{N+1}|^2\right).
\end{align*}
Next, utilizing Lemma \ref{lemm2} and \eqref{deltacap}, we obtain
\begin{equation*}
 \dot L_h(t) =  \dot E^{\rm ORFD}_h(t)+ \delta  \dot \phi_h(t)  \le-\delta (1-\delta L\eta)L_h(t) .
\end{equation*}
Finally, applying Gr\"onwall's inequality and Lemma \ref{lemm1}, we conclude
\begin{equation*}
(1-\delta L\eta) E^{\rm ORFD}_h(t) \le L_h(t)
\le (1+\delta L\eta)e^{-\delta (1-\delta L\eta) t} E^{\rm ORFD}_h(t).
\end{equation*}
Thus, \eqref{mainthm1} follows.
\end{proof}
Notice that the parameters $\sigma$, $\delta$, and $M$ in Theorem \ref{thm1} are identical to those obtained above for $\epsilon=1$. The maximum decay rate, given by $\sigma_{\max}(\delta) = \frac{1}{4\eta L}$, is achieved when $\delta(\xi_1, \xi_2) = \frac{1}{2\eta L}$. This result demonstrates the robustness of the model reduction \eqref{ORFD2} with respect to both the model reduction process and the discretization parameter $h$.

\subsection{Convergence of the Proposed Model Reduction to the PDE model as \texorpdfstring{$h\to 0$}{h->0}}
 We now establish the convergence of the semi-discretized model in \eqref{ORFD2} to the continuous PDE system \eqref{or1} as the discretization parameter $h$ tends to zero.

To achieve this, we define the error functions $\{\theta^1_j,\theta^2_j,\eta^1_j,\eta^2_j\}$, which measure the difference between the solutions of the discretized model and those of the PDE
\begin{equation*}
\begin{cases}
\theta^k_j(t) := {U}^k_j(t) - u^k(x_j,t), \quad \eta^k_j(t) := { W }^k_j(t) - w^k(x_j,t), & k=1,2,
\end{cases}
\end{equation*}
where ${ U^1_j, U^2_j, W^1_j,W^2_j}$ solve \eqref{discmodelwremainder}, while ${u^1,u^2,w^1,w^2}$ solve \eqref{ORFD2}.

The error functions ${\theta^1_j,\theta^2_j,\eta^1_j,\eta^2_j}$ satisfy the following system of equations
\begin{equation}
\label{errsstm}
\begin{cases}
\rho {\dot {\eta}}^{1}_{j+\frac{1}{2}}-\alpha\, \delta_x[\theta^1_{j+\frac{1}{2}}]+ \gamma \beta\; \delta_x [\theta^2_{j+\frac{1}{2}}]= s^1_{j+\frac{1}{2}}(t),&\\
\mu {\dot \eta}^{2}_{j+\frac{1}{2}}-\beta\,\delta_x [\theta^2_{j+\frac{1}{2}}] +\gamma \beta\,  \delta_x[\ \theta^1_{j+\frac{1}{2}}]=s^2_{j+\frac{1}{2}}(t),&\\
{\dot {\theta}}^{1}_{j+\frac{1}{2}} - \delta_x[ \eta^1_{j+\frac{1}{2}}]=r^1_{j+\frac{1}{2}}(t),
{\dot { \theta}}^{2}_{j+\frac{1}{2}}- \delta_x[ \eta^2_{j+\frac{1}{2}}]=r^2_{j+\frac{1}{2}}(t), & j=1,\ldots, N\\
\eta^1_0= \eta^2_0=0,\quad \alpha \theta^1_{N+1}-\gamma \beta  \theta^2_{N+1}= -k_1  \eta^1_{N+1},&\\
\beta  \theta^2_{N+1}-\gamma\beta  \theta^1_{N+1}=-k_2 \eta^2_{N+1}, &  t\in \mathbb{R}^+\\
( \theta^1_j,  \theta^2_j,  \eta^1_j,  \eta^2_j)(0) = \bm 0, &\hspace*{-1cm} j=0,1,\ldots, N+1.
\end{cases}
\end{equation}
  The discretized energy $e_h(t)$ for the solutions  $({\theta}^1, {\theta}^2, {\eta}^1,{ \eta}^2)$ of the error equation \eqref{errsstm}  is defined as follows
  \begin{equation}
\label{enerror}
	 e_h(t):=
	 \frac{h}{2} \sum\limits_{j=0}^N
\left(\rho     |\eta^1_{j+\frac{1}{2}}(t)|^2  +\mu    |\eta^2_{j+\frac{1}{2}}|^2 +\beta  |\gamma  \theta^1_{j+\frac{1}{2}}-\theta^2_{j+\frac{1}{2}}| ^2
+\alpha_1| \theta^1_{j+\frac{1}{2}}| ^2\right).
\end{equation}
Due to the positive-definite nature of the energy function,  when the energy of the error system \eqref{errsstm} converges to zero, the error itself must converge to zero as well.
\begin{lemma}
The energy $e_h(t)$ in \eqref{enerror} satisfies
\begin{equation*}
    \begin{split}
    \dot e_h(t)=&   -k_1\left({\eta} _{N+1}^1\right){}^2-k_2\left({ \eta} _{N+1}^2\right){}^2  +h \sum\limits _{j=0}^N \left\{ {\eta}_{j+\frac{1}{2}}^2  s_{j+\frac{1}{2}}^2 +{ \eta} _{j+\frac{1}{2}}^1  s_{j+\frac{1}{2}}^1\right.\\
    & +  r_{j+\frac{1}{2}}^1 \left(\alpha  {\theta} _{j+\frac{1}{2}}^1-\gamma \beta  {\theta} _{j+\frac{1}{2}}^2\right)\left.+  r_{j+\frac{1}{2}}^2 \left(\beta {\theta} _{j+\frac{1}{2}}^2-\gamma \beta {\theta} _{j+\frac{1}{2}}^1\right) \right\} .
    \end{split}
\end{equation*}
\end{lemma}
\begin{proof} Taking the time derivative of the discrete energy $e_h(t)$ along the solutions of \eqref{errsstm}, we obtain
\begin{equation*}
\begin{split}
    \dot e_h(t) =& h \sum\limits _{j=0}^N   \left\{\rho {\eta} _{j+\frac{1}{2}}^1 {\dot{\eta}}_{j+\frac{1}{2}}^{1} +  \mu    \eta _{j+\frac{1}{2}}^2 {\dot{\eta}}_{j+\frac{1}{2}}^{2}+  \alpha_1   {\theta} _{j+\frac{1}{2}}^1 {\dot{\theta}}_{j+\frac{1}{2}}^{1} \right.\\
&  \left. +  \beta  \left(\gamma  {\theta} _{j+\frac{1}{2}}^1-{\theta} _{j+\frac{1}{2}}^2\right)\left(\gamma  {\dot{\theta}} _{j+\frac{1}{2}}^{1}-{\dot{\theta}} _{j+\frac{1}{2}}^{2}\right) \right\}\\
      =&  h{\eta}^1 _{j+\frac{1}{2}} \left(\alpha\, \delta_x[\theta^1_{j+\frac{1}{2}}]- \gamma \beta\; \delta_x [\theta^2_{j+\frac{1}{2}}]+s^1_{j+\frac{1}{2}}(t)\right)\\
      &+h{\eta}^2 _{j+\frac{1}{2}} \left(\beta\,\delta_x [\theta^2_{j+\frac{1}{2}}] -\gamma \beta\,  \delta_x[\ \theta^1_{j+\frac{1}{2}}]+s^2_{j+\frac{1}{2}}(t)\right)\\
      & +\sum\limits _{j=0}^N  \alpha_1 h {\theta}^1_{j+\frac{1}{2}} \left( \delta_x[ \eta^1_{j+\frac{1}{2}}]+r^1_{j+\frac{1}{2}}(t)\right)\\
      & \hspace*{-1cm} +  \sum\limits _{j=0}^N \beta h \left(\gamma  {\theta} _{j+\frac{1}{2}}^1-{\theta} _{j+\frac{1}{2}}^2\right)\left[\gamma\left( \delta_x[ \eta^1_{j+\frac{1}{2}}]+r^1_{j+\frac{1}{2}}(t)\right)\right. -\left. \left( \delta_x[ \eta^2_{j+\frac{1}{2}}]+r^2_{j+\frac{1}{2}}(t)\right) \right].
\end{split}
\end{equation*}
By Lemma \ref{summa}, we have the following identities
\begin{equation*}
\begin{cases}
 \alpha h  \sum\limits_{j=0}^N \left( \eta _{j+\frac{1}{2}}^1 \delta _x\left[\theta ^1{}_{j+\frac{1}{2}}\right] + \theta _{j+\frac{1}{2}}^1 \delta _x\left[\eta _{j+\frac{1}{2}}^1\right]\right)=\alpha  \left(\eta _{N+1}^1 \theta _{N+1}^1-\eta _0^1 \theta _0^1\right),\\
  \beta h\sum\limits_{j=0}^N  \left(\eta _{j+\frac{1}{2}}^2 \delta _x\left[\theta ^2{}_{j+\frac{1}{2}}\right] + \theta _{j+\frac{1}{2}}^2 \delta _x\left[\eta _{j+\frac{1}{2}}^2\right]\right)=\beta  \left(\eta _{N+1}^2 \theta _{N+1}^2-\eta _0^2 \theta _0^2\right),\\
  - \beta\gamma h \sum\limits_{j=0}^N \left( \eta _{j+\frac{1}{2}}^1 \delta _x\left[\theta ^2{}_{j+\frac{1}{2}}\right] + \theta _{j+\frac{1}{2}}^2 \delta _x\left[\eta _{j+\frac{1}{2}}^1\right]\right) =-\beta\gamma  \left(\eta _{N+1}^1 \theta _{N+1}^2+\eta _0^1 \theta _0^2\right),\\
      -\beta\gamma h\sum\limits_{j=0}^N \left( \eta _{j+\frac{1}{2}}^2 \delta _x\left[\theta ^1{}_{j+\frac{1}{2}}\right] + \theta _{j+\frac{1}{2}}^1 \delta _x\left[\eta _{j+\frac{1}{2}}^2\right]\right) =-\beta\gamma  \left(\eta _{N+1}^2 \theta _{N+1}^1+\eta _0^2 \theta _0^1\right).
\end{cases}
\end{equation*}
By substituting these in $  \dot e_h(t)$ above leads to
\begin{equation*}
\begin{split}
 \dot e_h(t) =&\eta _{N+1}^1 \left(\alpha \theta _{N+1}^1-\gamma \beta  \theta _{N+1}^1\right)+\eta _{N+1}^2 \left(\beta \theta _{N+1}^2-\gamma \beta \theta _{N+1}^2\right)\\
 &+h \sum\limits _{j=0}^N \left\{ {\eta}_{j+\frac{1}{2}}^2  s_{j+\frac{1}{2}}^2  +{ \eta} _{j+\frac{1}{2}}^1 s_{j+\frac{1}{2}}^1 +r_{j+\frac{1}{2}}^1 \left(\alpha  {\theta} _{j+\frac{1}{2}}^1-\gamma \beta  {\theta} _{j+\frac{1}{2}}^2\right)\right.\\
  &   \left.+ r_{j+\frac{1}{2}}^2 \left(\beta {\theta} _{j+\frac{1}{2}}^2-\gamma \beta {\theta} _{j+\frac{1}{2}}^1\right)\right\}.
\end{split}
\end{equation*}
Finally, applying the boundary conditions from \eqref{errsstm} completes the proof.
\end{proof}

  \subsection{Exponential Stability of the Error Dynamics by the Lyapunov Approach}
 For $\tilde\delta>0,$ define the Lyapunov functional as the following
  \begin{equation} \label{psidef}
  l_h(t)=e_h(t)+ \tilde\delta \psi_h(t):=e_h(t)+  h\sum\limits_{j=0}^N   \rho  x_{j+\frac{1}{2}}
        \theta_{j+\frac{1}{2}}^{1}
        \eta_{j+\frac{1}{2}}^{1}
        +
        \mu x_{j+\frac{1}{2}}
        \theta_{j+\frac{1}{2}}^{2}.
\end{equation}
Analogous to Lemma \ref{lemm2}, for $\tilde \delta<\frac{1}{L\eta}$, we obtain the following bound
 \begin{equation}\label{lemm111} (1-\tilde\delta L \eta) e_h(t) \le l_h(t) \le (1+\tilde\delta L \eta) e_h(t). \end{equation} The inequality follows from standard Lyapunov energy estimates and the structure of $\psi_h(t)$.

\begin{lemma}\label{lemm12}$\psi(t)$ in \eqref{psidef} satisfies the following estimate
\begin{align*}
 \dot \psi_h(t)\leq&  - e_{h}(t)+
    \frac{1}{2} \left\{
    \left(
     \rho+   \frac{k_{1}^{2}}{\alpha_1} \right) |\eta^1_{N+1}|^2 + \left(\mu + \frac{k_{2}^2(\alpha + \gamma^2\beta) }{2\alpha_1 \beta}  \right) |\eta^2_{N+1}|^2\right\}\\
    &+ h \sum\limits_{j=0}^N\x \left(\rho  \veI \rI +
   \mu \veII \rII\right) \\
   &+h \sum\limits_{j=0}^N \x\left(  \ueI \sI+   \ueII \sII \right).
\end{align*}
\end{lemma}

\begin{proof}
By a direct calculation with \eqref{errsstm}
\begin{align*}
   \dot \psi_h(t)= & h \sum\limits_{j=0}^N \x \left\{    \rho   \delta_x
            [ \eta^{1}]
            \eta^{1}+    \rho   \eta^{1}_{j+\frac{1}{2}}
            r^{1}_{j+\frac{1}{2}} \right. +  \alpha \delta_x [ \theta^{1}_{j+\frac{1}{2}}]  \theta^{1}_{j+\frac{1}{2}} \\
            &  -
            \gamma
            \beta  \theta^{1}_{j+\frac{1}{2}}
            \delta_x
            [\theta^{2}_{j+\frac{1}{2}}] +  \theta^{1}_{j+\frac{1}{2}} s^{1}_{j+\frac{1}{2}}+
            \mu
            \delta_{x}
            [
                \eta^{2}_{j+\frac{1}{2}}
            ]
            \eta^{2}_{j+\frac{1}{2}}\\
            & +  \mu   \eta^{2}_{j+\frac{1}{2}}r^2_{j+\frac{1}{2}} + \beta  \delta_{x}[\theta^{2}_{j+\frac{1}{2}}   ]\theta^{2}_{j+\frac{1}{2}}h \beta\gamma \delta_{x} [ \theta^{1}_{j+\frac{1}{2}}] \theta^{2}_{j+\frac{1}{2}}
       \left. +
                   \theta^{2}_{j+\frac{1}{2}}s^{2}_{j+\frac{1}{2}}\right\}.
\end{align*}
Next, we use the summation by parts formula in Lemma \ref{summa}
\begin{align*}
\dot \psi_h(t)
    =  &  -2e_h(t) + h \sum\limits_{j=0}^N \x\left\{
            -\rho
            \delta_{x}  [ \eta^{1}_{j+\frac{1}{2}}  ]
            \eta^{1}_{j+\frac{1}{2}}  \right.  \\
            &-
            \mu
            \delta_{x}
            [
                \eta^{2}_{j+\frac{1}{2}}
            ]
            \eta^{2}_{j+\frac{1}{2}}
        -
            \alpha
            \delta_{x}
            [
                \theta^{1}_{j+\frac{1}{2}}
            ]
            \theta^{1}_{j+\frac{1}{2}}
 -
            \beta
            \delta_{x}
            [
                \theta^{2}_{j+\frac{1}{2}}
            ]
            \theta^{2}_{j+\frac{1}{2}} + \beta \gamma
            \delta_{x}
            [
                \ueI
            ]
            \ueII\\
            & \left.+
            \beta
            \gamma
                    \delta_{x}
            [
                \ueII
            ]
            \ueI +
            \rho
            \veI
            \rI
        +  \mu
            \veII
            \rII      +
            \ueI
            \sI
        +
            \ueII
            \sII  \right\}\\
            &-
        \frac{h}{4}
        \sum\limits_{j=0}^N
           \left\{ \rho
           |
                \eta^{1}_{j+1}-\eta^{1}_{j}
                |^2
        -
            \mu
            |
                \eta^{2}_{j+1}-\eta^{2}_{j}
                |^2 \right.
    \\
    &-
            \alpha
            |
                \theta^{1}_{j+1}-\theta^{1}_{j}
                |^2
        -
                 \beta
                 |
                \theta^{2}_{j+1}-\theta^{2}_{j}
                |^2
+
       2 \beta
            \gamma
            (
                \theta^{1}_{j+1}-\theta^{1}_{j}
            )
            (
                \theta^{2}_{j+1}-\theta^{2}_{j}
            ) \\
            & +
        \alpha
        |\theta^{1}_{N+1}|^2
        +
        \beta
        |\theta^{2}_{N+1}|^2
        -
        2
        \beta
        \gamma
        \theta^{1}_{N+1}
        \theta^{2}_{N+1} \left.
        +
        \rho
        |\eta^1_{j+\frac{1}{2}}|^2
        +
        \mu
        |\eta^2_{j+\frac{1}{2}}|^2 \right\}.
    \end{align*}
    Therefore,
    \begin{align*}
 2\dot \psi_h(t)=  &  -  2e_{h}(t)
        \\
        &+
        2   h \sum\limits_{j=0}^N
                \x
  \left( \rho
                             \veI
                \rI+
                \mu
                \veII
                \rII  +
                \ueI
                \sI
            +  \ueII
                \sII
        \right)
    \\ \nonumber
    &
        +\frac{h}{4}
        \sum\limits_{j=0}^N
            \left\{-\rho
            |
                \eta^{1}_{j+1}-\eta^{1}_{j}
                |^2
        - \mu
        |
                \eta^{2}_{j+1}-\eta^{2}_{j}
                |^2\right.
    \\ \nonumber
    &-    \alpha
    |
                \theta^{1}_{j+1}-\theta^{1}_{j}
                |^2
        -\beta
        |
                \theta^{2}_{j+1}-\theta^{2}_{j}
                |^2
  \left. +
       2 \beta
                \gamma
                (
                    \theta^{1}_{j+1}-\theta^{1}_{j}
                )
                (
                    \theta^{2}_{j+1}-\theta^{2}_{j}
                )\right\}
    \\ \nonumber
    &+
        \alpha
        |\theta^{1}_{N+1}|^2
        +
        \beta
        |\theta^{2}_{N+1}|^2
        -
        2
        \beta
        \gamma
        \theta^{1}_{N+1}
        \theta^{2}_{N+1}
        +
        \rho
        |\eta^{1}_{N+1}|^2
        +
        \mu
        |\eta^{2}_{N+1} |^2 .
\end{align*}
Simplifying the equations above and dropping the negative terms yield
\begin{align*}
\dot \psi_h(t)\leq&
    - e_{h}(t)
    + h \sum\limits_{j=0}^N \x\left\{ \rho  \veI \rI + \mu \veII \rII+ \ueI \sI    +  \ueII \sII     \right\}
\\ \nonumber
& +
    \frac{1}{2}
    \left[
        \alpha |\theta^{1}_{N+1}|^2
        - 2 \beta \gamma \theta^{1}_{N+1}\theta^{2}_{N+1}
        + \beta |\theta^{2}_{N+1}|^2
        + \rho |\eta^{1}_{N+1}|^2
        + \mu |\eta^{2}_{N+1} |^2
    \right].
\end{align*}
Finally, the conclusion follows by adopting the boundary conditions in \eqref{errsstm}.
\end{proof}

\begin{theorem} \label{thm111}[Convergence of Solutions] Define the positive constants
\begin{align*}
 C_1:=& \max \left(
     \frac{1}{\rho} + \epsilon L \sqrt{\frac{\rho}{\alpha_1}},  \frac{1}{\mu} + \epsilon L\sqrt{\frac{\mu}{\beta}}\left(1+ \gamma\frac{\beta}{\alpha_1}\right),  1 + \epsilon L \left(\sqrt{ \frac{\rho}{ \alpha_1} } +\gamma\sqrt{ \frac{ \mu}{\alpha_1} }  \right), \right. \\
     & 1 + \epsilon L \gamma\sqrt{ \frac{  \mu}{ \beta} }, 1 + \epsilon L \left(\gamma\sqrt{ \frac{1}{\rho \alpha_1} } +\sqrt{ \frac{ 1}{\mu \alpha_1} }  \right),\left.     \frac{1}{\mu} + \epsilon L\left(\frac{\gamma}{\sqrt{\alpha_1\mu}}+\frac{1}{\sqrt{\beta \mu}} \right)
    \right),\\
 C_2:=&  \max\left(\rho, \mu, \alpha +\gamma^2\beta, 2\beta\right),
\end{align*}
and choose
\begin{equation}\label{tildedeltacap}
   \tilde \delta<\min \left(\frac{1}{L\eta},\frac{2k_1\alpha_1}{\alpha_1\rho + k_1^2}, \frac{4k_2\beta}{2\alpha_1\beta\mu+(1+\frac{2\gamma^2 \beta}{\alpha_1})k_2^2}\right).
\end{equation}
If $u^1(x,t),u^2(x,t),w^1(x,t),w^2(x,t)\in C^1([0,\infty);C^3[0,1])$ solve the system \eqref{or1} and $\left\{(u^1, u^2, w^1, w^2)_j(t)\right\}_{j=0}^{N+1}$  solves the system \eqref{ORFD2}, then the error functions \linebreak $\left\{(\theta^1, \theta^2, \eta^1, \eta^2)_j(t)\right\}_{j=0}^{N+1}$  satisfy the system \eqref{errsstm} and the following  exponential stability result holds
\begin{equation}\label{expstab}
  e_h(t) \;\leq \frac{1+\tilde\delta L\eta}{1-\tilde\delta L\eta}e^{ -\frac{\tilde\delta (1-\tilde \delta L\eta )}{2}t} e_h(0)+   \frac{2C_1^2 C_2 D^2 h^4}{\tilde\delta^2(1-\tilde\delta L\eta)^2},
\end{equation}
where the positive constant $D$ is defined along the proof.
\end{theorem}

\begin{proof}
Since $u^1(x,t),u^2(x,t),v^1(x,t),v^2(x,t)\in C^1([0,\infty);C^3[0,1]),$ by \eqref{remainders}, \linebreak there exists a constant $D>0$, independent of $h$, such that
\begin{equation*}
    |\rI|,  |\rII|,  |\sI| , |\sII|\leq Dh^2
\end{equation*}
for all $ j= 0,1,...,N+1$, and $t\in [0,\infty]$.
By Lemma \ref{lemm12},
\begin{align*}
 \dot l_h(t)
\leq& -\tilde\delta e_{h}(t)
    -
    \left(
        k_{1} - \frac{\tilde \delta}{2}\left(\rho +\frac{k_{1}^{2}}{\alpha_1}\right)
    \right)
    |\eta^1_{N+1}|^2 \\
   & -
    \left(
    k_{2}
    - \frac{\tilde\delta}{2}
    \left( \mu +\frac{k_{2}^2(\alpha_1 +\beta \gamma^2) }{2\alpha_1 \beta} \right)\right)|\eta_{N+1}^2|^2
\\ \nonumber
&+
    h \sum\limits_{j=0}^N \left\{ \veII \sII + \veI \sI + \rI ( \alpha \ueI - \gamma \beta \veII )
    \right.\\
    & \hspace*{3cm}\left.+
    \rII (\beta \ueII - \gamma \beta \ueI ) \right\}\\
    &+ \tilde\delta h \sum\limits_{j=0}^N \x \left\{\rho  \veI \rI
    + \mu \veII \rII +
   \ueI \sI
   +\ueII \sII \right\}.
\end{align*}
By \eqref{tildedeltacap},
\begin{align*}
\dot l_h(t)=&
    -\tilde\delta e_{h}(t)+
    h  \sum\limits_{j=0}^N \left\{ \veII \sII
    + \veI \sI
    + \alpha_1 \ueI \rI  \right.
\\ \nonumber
&  \left.+  \beta( \gamma \rI -\rII) (\gamma \ueI - \ueII )\right\}
 + \tilde\delta h L \sum\limits_{j=0}^N \left\{ \rho \veI \rI +\ueI \sI \right.
\\ \nonumber
&  + \mu \veII \left(\rII -\gamma\rI +\gamma\rI\right)  \left.+ \left(\ueII -\gamma \ueI +\gamma \ueI \right) \sII. \right\}
\end{align*}
By several algebraic manipulations and  the generalized Young's inequality for $\xi>0$,
\begin{align*}
 \dot l_h(t)=&
    -\tilde\delta e_{h}(t)
+    \frac{h}{2}\left[
    \frac{1}{\xi}
 \left\{  \sum\limits_{j=0}^N \left(
     \frac{1}{\rho} + \tilde\delta L \sqrt{\frac{\rho}{\alpha_1}}
    \right)
     \rho|\veI|^2 \right.\right.
\\ \nonumber
&+
    \left(
     \frac{1}{\mu} + \tilde\delta L\sqrt{\frac{\mu}{\beta}}\left(1+ \gamma\frac{\beta}{\alpha_1}\right)
    \right)
\mu |\veII|^2
+
  \left(
1 + \tilde\delta L \gamma\sqrt{ \frac{  \mu}{ \beta} }    \right)
\beta |\gamma \ueI - \ueII |^2
\\ \nonumber
&\left. +
     \left(
1 + \tilde\delta L \left(\gamma\sqrt{ \frac{1}{\rho \alpha_1} } +\sqrt{ \frac{ 1}{\mu \alpha_1} }  \right)  \right)    \frac{h}{2}
\alpha_1 |\ueI|^2 \right\}\\
& +
    \xi  \left\{ \sum\limits_{j=0}^N
    \left(
      \frac{1}{\rho}  + \tilde\delta L\frac{1}{\sqrt{\alpha_1\rho}} \right)
  \rho|\sI|^2
    +
    \left(
          \frac{1}{\mu} + \tilde\delta L\left(\frac{\gamma}{\sqrt{\alpha_1\mu}}+\frac{1}{\sqrt{\beta \mu}} \right)
    \right)
 \mu|\sII|^2
 \right.\\
&+
    \left(
1 + \tilde\delta L \gamma\sqrt{ \frac{  \mu}{ \beta} }
    \right)
 \beta | \gamma \rI -\rII|^2 \\
 &\left.  \left.+
    \left(
      1 + \tilde\delta L \left(\sqrt{ \frac{\rho}{ \alpha_1} } +\gamma\sqrt{ \frac{ \mu}{\alpha_1} }  \right)
    \right)
\alpha_1 |\rI|^2 \right\} \right].
\end{align*}
Now choose $\xi= \frac{2C_1}{\tilde\delta}$ so that
\begin{align*}
 \dot l_h(t)
 \le  &    -\frac{\tilde\delta (1-\tilde \delta L\eta )}{2}l_{h}(t) \\
 & +
    \frac{C_1^2 h}{\tilde\delta}
    \left(
    \sum\limits_{j=0}^N\left\{ \rho|\sI|^2
    \right.\right. \left.  +
    \mu|\sII|^2
    +
 \beta | \gamma \rI -\rII|^2
    +
 \alpha_1 |\rI|^2
    \right)
\end{align*}
where \eqref{lemm111} is used. By Gr\"onwall's inequality,
\begin{align*}
l_h(t)\le& e^{ -\frac{\tilde\delta (1-\tilde \delta L\eta )}{2}t} l_h(0)+   \frac{C_1^2 h}{\tilde\delta}
 \int_0^t e^{ -\frac{\tilde\delta (1-\tilde \delta L\eta )}{2}(t-\tau)} \left(
    \sum\limits_{j=0}^N\left\{ \rho|\sI|^2
  \right.\right.\\
  & \left.  +
    \mu|\sII|^2
    +
 \beta | \gamma \rI -\rII|^2
    +
 \alpha_1 |\rI|^2
    \right)
\end{align*}
By \eqref{lemm111} one more time,
\begin{align*}
e_h(t)\le &\frac{1+\tilde\delta L\eta}{1-\tilde\delta L\eta}e^{ -\frac{\tilde\delta (1-\tilde \delta L\eta )}{2}t} e_h(0) \\
&+   \frac{C_1^2 D^2 h^4 \max\left(\rho, \mu, \alpha +\gamma^2\beta, 2\beta\right)}{\tilde\delta(1-\tilde\delta L\eta)}
 \int_0^t e^{ -\frac{\tilde\delta (1-\tilde \delta L\eta )}{2}(t-\tau)} d\tau\\
 \le &\frac{1+\tilde\delta L\eta}{1-\tilde\delta L\eta}e^{ -\frac{\tilde\delta (1-\tilde \delta L\eta )}{2}t} e_h(0)+   \frac{2C_1^2 C_2 D^2  h^4}{\tilde\delta^2(1-\tilde\delta L\eta)^2}.
\end{align*}
This makes the end of the proof.
\end{proof}

\subsection{Convergence of Discretized Energy to the Energy of the PDE  as \texorpdfstring{$h\to 0$}{h->0}}
In the prior section, the convergence of the semi-discretized model is proved as $h\to 0.$ Now, we are ready to show the convergence of this discretized energy $E_h^{\rm ORFD}(t)$ in \eqref{eq4} to the energy $E(t)$ in  \eqref{enpde2} of the PDE model.

\begin{theorem}[Convergence of the Energy]
Suppose that $(u^1, u^2, w^1, w^2)(x,t)\in \left[C^1\left([0,\infty);C^3[0,1]\right)\right]^4$ solve the PDE system~\eqref{or1}, and let $\left\{(u^1, u^2, w^1, w^2)_j(t)\right\}_{j=0}^{N+1}$ solve the discretized system~\eqref{ORFD2}. Then, the corresponding error functions\\ $\left\{(\theta^1, \theta^2, \eta^1, \eta^2)_j(t)\right\}_{j=0}^{N+1}$ satisfy the error system~\eqref{errsstm}.

Assume that the initial conditions in the discretized model~\eqref{ORFD2} are exact and that the initial discretized energy $E_h^{\rm ORFD}(0)$ defined by~\eqref{eq4} converges to the energy $E(0)$ defined by~\eqref{enpde2} of the PDE model as $h\to 0$. Then, the following energy convergence estimate holds:
\begin{equation}\label{con4}
\left|E_h^{\rm ORFD}(t)-E(t)\right| \le \left|E_h^{\rm ORFD}(0)-E(0)\right|+ c(k_1,k_2)\,h^{3/2}\sqrt{t}, \quad t\ge 0,
\end{equation}
where the positive constant $c(k_1,k_2)$ (explicitly defined in the proof) is independent of $h$. Consequently, the energy convergence in the uniform norm follows as
\begin{equation}\label{convergence}
\lim_{h\to 0}\|E_h^{\rm ORFD}(t)-E(t)\|_{C[0,\infty]}=0.
\end{equation}
\end{theorem}
\begin{proof}
Since the initial conditions are exact, then $e_h(0)\equiv 0,$ and therefore, by \eqref{expstab}, there exists a constant $c_1$ such that
\begin{equation}\label{con1}
\begin{split}
&\frac{h}{2} \left(\left|w^1_{N+1}(t)-w^1(L,t)\right|^2+ \left| w^2_{N+1}-w^2(L,t)\right|^2 \right)\\
&\qquad \le \frac{h}{2} \left(\left|\eta^1_{N+1}(t)\right|^2+ |\eta^2_{N+1}(t)|^2\right)\le 4\min\left(\frac{1}{\rho},\frac{1}{\mu}\right)e_h(t)\le c_1 h^4.
    \end{split}
\end{equation}
By Theorem \ref{thm1} and Lemma \ref{lemm1}, respectively,
\begin{align*}
E(t)&=E(0) -\int_0^t\left(k_1\left|w^1(L,t)\right|^2 + k_2 \left|w^2(L,t)\right|^2 \right) dt,\\
E_h^{\rm ORFD}(t)&=E_h^{\rm ORFD}(0) -\int_0^t \left(k_1\left|w^1_{N+1}(t)\right|^2 + k_2 \left|w^2_{N+1}(t)\right|^2\right) dt,
\end{align*}
and since  $E_h(t), E_h^{\rm ORFD}(t)\to 0$ exponentially as $h\to 0,$
\begin{align*}
E(0) =&\int_0^t\left(k_1\left|w^1(L,t)\right|^2 + k_2 \left|w^2(L,t)\right|^2 \right)dt,\\
E_h^{\rm ORFD}(0)=& \int_0^t\left(k_1\left|w^1_{N+1}(t)\right|^2 + k_2 \left|w^2_{N+1}(t)\right|^2 \right)dt.
\end{align*}
Since it is assumed that $E_h^{\rm ORFD}(0) \to E(0)$ as $h\to 0, $ both integrals terms above, are bounded, i.e. there exits $c_2,c_3>0$ such that
\begin{equation}\label{con2}
\int_0^t(|w^1(L,t)|^2 +  |w^1_{N+1}(t)|^2 )dt<c_2,\quad \int_0^t(|w^2(L,t)|^2 + |w^2_{N+1}(t)|^2 )dt<c_3.
\end{equation}
Therefore, for any $t\ge 0, $ by the H\"older's and Cauchy-Scwartz's inequalities,
\begin{align*}
&\left|E_h^{\rm ORFD}(t)-E(t)\right|\le \left|E_h^{\rm ORFD}(0)-E(0)\right|\\
&~\quad + k_1\left|\int_0^t \left(\left|w^1_{N+1}(t)\right|^2-\left|w^1(L,t)\right|^2\right)dt\right|
+k_2 \left|\int_0^t \left(\left|w^2_{N+1}(t)\right|^2-\left|w^2(L,t)\right|^2\right)dt \right|\\
&\quad\le \left|E_h^{\rm ORFD}(0)-E(0)\right|\\
&\qquad + k_1 \left[\int_0^t \left|w^1_{N+1}(t)-w^1(L,t)\right|^2 dt \right]^{\frac{1}{2}}
\left[\int_0^t \left|w^1_{N+1}(t)+w^1(L,t)\right|^2 dt \right]^{\frac{1}{2}}\\
&\qquad + k_2\left[\int_0^t \left|w^2_{N+1}(t)-w^2(L,t)\right|^2 dt\right]^{\frac{1}{2}}
\left[\int_0^t \left|w^2_{N+1}(t)+w^2(L,t)\right|^2 dt\right]^{\frac{1}{2}}\\
&\quad\le \left|E_h^{\rm ORFD}(0)-E(0)\right|\\
&\qquad + \sqrt{2}k_1 \left[\int_0^t \left|w^1_{N+1}(t)-w^1(L,t)\right|^2 dt\right]^{\frac{1}{2}}
\left(\int_0^t \left[\left|w^1_{N+1}(t)\right|^2+\left|w^1(L,t)\right|^2\right) dt\right]^{\frac{1}{2}}\\
&\qquad + \sqrt{2}k_2\left[\int_0^t \left|w^2_{N+1}(t)-w^2(L,t)\right|^2 dt\right]^{\frac{1}{2}}
\left(\int_0^t \left[\left|w^2_{N+1}(t)\right|^2+\left|w^2(L,t)\right|^2\right) dt\right]^{\frac{1}{2}}.
\end{align*}
Utilizing \eqref{con1} and \eqref{con2} lead to
\begin{align*}
\left|E_h^{\rm ORFD}(t)-E(t)\right| \le \left|E_h^{\rm ORFD}(0)-E(0)\right| +\sqrt{2c_1} \left(k_1\sqrt{c_2} + k_2\sqrt{c_3} \right) h^{3/2} \sqrt{t}.
\end{align*}
By defining the constant $c(k_1,k_2):=\sqrt{2c_1} \left(k_1\sqrt{c_2} + k_2\sqrt{c_3} \right),$  \eqref{con4} follows.

Finally, since both energies $E(t)$ and $E_h^{\rm ORFD}(t)$ exhibit exponential decay (Theorem \ref{thm1} and Theorem \ref{mainexp}), for sufficiently large $t'>t$, they become sufficiently small, thereby establishing the desired convergence result.
\end{proof}

\section{Simulations and Numerical Experiments}
\label{Simu}

In this section, numerical simulations and experiments are conducted using Wolfram Mathematica to demonstrate the effectiveness of the proposed finite-difference-based model reduction~\eqref{ORFD2}, introduced in Section~\ref{ModelRed}, compared to the finite element-based model reduction introduced in Section~\ref{FEMSplines}. The simulations utilize realistic piezoelectric material parameters as provided in Table~\ref{tablo}. It is important to highlight that employing realistic material parameters as presented in Table~\ref{tablo} introduces practical complexities. Specifically, the significant contrast between mechanical and electromagnetic properties, characterized by relations such as $\mu \ll \rho$ and $\frac{\alpha}{\rho}\ll \frac{\beta}{\mu}$, may compromise numerical robustness and stability. Additionally, the speed $\sqrt{\frac{\beta}{\mu}}$ closely approximates the speed of light, further emphasizing these inherent numerical challenges.

\begin{table}[htb!]
  \caption{Realistic Piezoelectric Material Constants}
  \centering
  \label{tablo}
  \begin{tabular}{l|c|c|c}
    \hline
    Property&Symbol&Value&Unit\\
    \hline
    Length of the beam&$L$&$1$&m \\
    Mass density&$\rho$&$6000$&kg/${\rm m}^3$\\
    Magnetic permeability&$\mu$&$10^{-6}$ & H/m\\
    Elastic stiffness&$\alpha$&$10^9$& N/${\rm m}^2$\\
    Piezoelectric constant&	$\gamma$&$10^{-3}$&C/${\rm m}^3$\\
    Impermittivity  &$\beta$& $ 10^{12}$& m/F\\
    \hline
  \end{tabular}
\end{table}

The selection of appropriate feedback amplifiers $k_1$ and $k_2$ is crucial for achieving rapid exponential stability. For the realistic material parameters provided in Table~\ref{tablo}, the maximal decay rate established in \cite[Theorem 1]{ozer2024exponential} is $\sigma_{\text{max}}=102.04$. We select feedback amplifiers within the proposed safe intervals as follows
\begin{equation}
k_1=10^6\in \left( 7.17 \times 10^5, 4.18 \times 10^6\right), ~~
k_2=10^6\in \left(1.02 \times 10^{-4}, 9.78 \times 10^9\right).
\end{equation}

Let $\left\{\lambda^{\pm}_{1j}, \lambda^{\pm}_{2j}\right\}(h)$ and $\left\{\bm{\Psi}^{\pm }_{1j}, {\bm{\Psi}^{\pm }_{2j}}\right\} (h)$ be the eigenvalues and eigenvectors corresponding to the first order formulation of the system \eqref{FEM}. In particular, they satisfy following eigenvalue problem that is $\eqref{disc}$ with feedback controllers,
\begin{equation}
\label{eq:controlEigproblem}
\begin{bmatrix}
  0& \bm I\\
  \mathbb{A}& \mathbb{K}
\end{bmatrix}\bm{\Psi}= \lambda \bm{\Psi}, \quad \text{where} \quad \begin{cases}
  \mathbb{A}= -\left(\bm C_1^{-1}\bm C_2\right)  \otimes \left(\bm M^{-1}\bm A_h\right),\\
  \mathbb{K}= -\left(\bm C_1^{-1}\bm C_3\right) \otimes \left(\bm M^{-1}\bm B\right),
\end{cases}
\end{equation}
 where the chosen mass matrix, $\bm{M}$ (defined in \eqref{eq:massFEM} for FEM and \eqref{eq:massOR} for ORFD), determines which approximation method is used.

\begin{remark}
  \label{rmk:filtering}
The direct Fourier filtering procedure described in \cref{sec:filtering} must be separately applied to each signed half of the eigenvalue branches, as outlined in \eqref{eq:filtered}. Although these branches are explicitly identifiable in the absence of control, distinguishing them in the controlled case is nontrivial. Therefore, we develop an algorithm based on the observation that introducing viscous damping in one wave equation significantly perturbs only one eigenvalue branch. The auxiliary eigenvalue problem
\begin{equation}
    \label{eq:testEigproblem}
    \begin{bmatrix}
      0& \bm I\\
      \mathbb{A} & \mathbb{K}_\epsilon
    \end{bmatrix}\hat{\bm \Psi}= \hat \lambda \hat{\bm \Psi}, \quad \text{where} \quad
    \mathbb{K}_\epsilon= \mathbb{K}-\left(\bm C_1^{-1}\begin{bmatrix}
        \epsilon& 0\\
        0& 0
      \end{bmatrix}\right) \otimes \left(\bm M^{-1}\right).
  \end{equation}
   incorporates these perturbations explicitly. The detailed steps of this branch separation algorithm are provided in \cref{alg:BranchSeperation}. Practically, eigenvalue comparisons based on imaginary components yield a robust and computationally efficient implementation of this method.

	\begin{algorithm}
		\caption{Separating branches of  eigenvalues and their eigenfunctions}
		\label{alg:BranchSeperation}
		\begin{algorithmic}[1]
		\State{Compute the actual eigenpairs $\left\{\lambda, \bm \Psi \right\}$, and test eigenpairs  $\{\hat \lambda, \hat{\bm \Psi} \}$from solving \eqref{eq:controlEigproblem} and \eqref{eq:testEigproblem}, respectively.}
    \For{$i = 1$ \textbf{to} $4N+4$}
      \State{$j= \operatorname{argmin}_j\norm{{\bm \Psi}_i -\hat{\bm \Psi}_j } $}
      \If{ $\left|\text{Re}(\lambda_i)-\text{Re}(\hat \lambda_j)\right|<\text{tol}_\epsilon$}
        \State \text{Assign the eigenpair $\{\lambda_i, {\bm \Psi}_i\}$ to the first branch}
      \Else
        \State \text{Assign the eigenpair $\{\lambda_i, {\bm \Psi}_i\}$ to the second branch}
      \EndIf
    \EndFor
	\end{algorithmic}
	\end{algorithm}
	\end{remark}

\begin{figure}[htb!]
	\graphicspath{{visuals/}}
	\includegraphics[width=0.48\columnwidth]{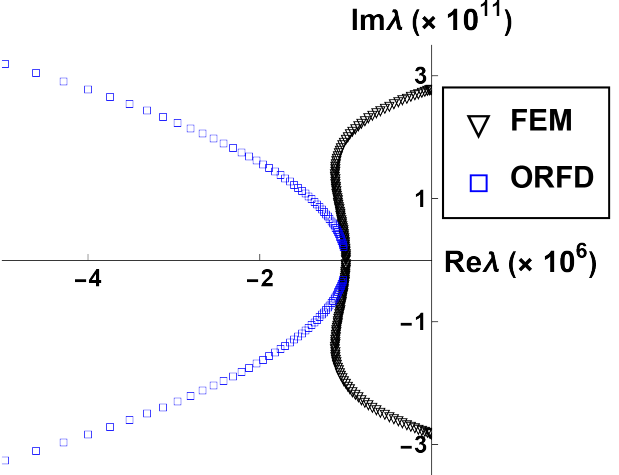} \quad \includegraphics[width=0.48\columnwidth]{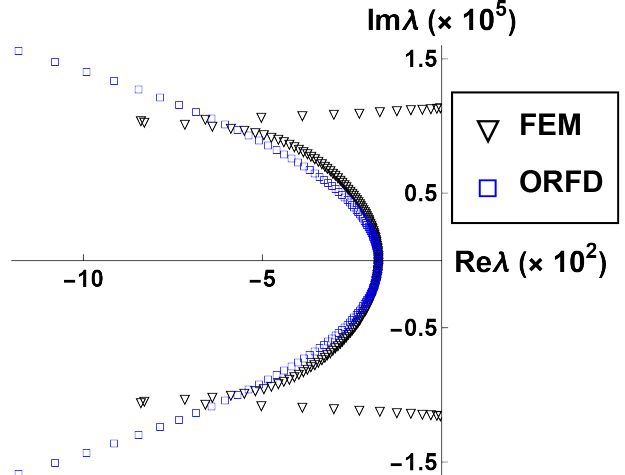}
	\caption{Spectral plots comparing FEM (black triangles) and ORFD (blue squares) eigenvalues for $N=80$. The ORFD eigenvalues remain robustly bounded away from the imaginary axis, unlike the FEM eigenvalues, which approach it and require numerical filtering.}
	\label{eigvalsmerged}
\end{figure}

With $N=80$, the spectral plots in Figure \ref{eigvalsmerged} illustrate the eigenvalues corresponding to the first-order formulation of \eqref{FEM} and \eqref{ORFD2}. Clearly, the Finite Element Method (FEM) eigenvalues approach the imaginary axis, whereas the Order-Reduction-Based Finite Differences (ORFD) eigenvalues remain uniformly bounded away from it. Due to the FEM model's high-frequency eigenvalues approaching the imaginary axis, numerical filtering is applied to ensure that the eigenvalues on each branch remain sufficiently bounded away from the imaginary axis.

A reliable measure of system stability is the maximum real part of the eigenvalues. We examine three different node counts: $N = 40, ~80, ~160$. As shown in \cref{tablo2}, in the FEM case without filtering, i.e., (i) $j^*=0$, the maximum real part of the eigenvalues approaches zero as $N$ increases, indicating severe instability in FEM solutions. With some filtering, i.e., (ii) $j^*=5$ and (iii) $j^*=10$, the maximum real parts of the eigenvalues shift away from the imaginary axis.  Slow exponential decay rates can also be observed.

\begin{table}[!htb]
	\label{tablo2}
	\centering
	\begin{tabular}{c|c|c|c|c}
		\toprule
		{\textbf{Method}}  & \textbf{Fourier Filtering} &\multicolumn{3}{c}{\textbf{Maximum real part of eigenvalues}} \\
		\midrule
		& $j^*$ &  N = 40 & N = 80 & N = 160\\
		\midrule
		&  0 &$-13.1571$&$-3.37871$&$-0.8557$\\
		FEM  & 	5 & $-177.033$&$-125.217$&$-30.9887$\\
		& 10 &$-177.033$&$-177.001$&$-106.099$\\
		\midrule
		ORFD& N/A  &$-177.055$&$-176.935$&$-174.897$\\
		\bottomrule
	\end{tabular}
		\caption{Maximum real part of eigenvalues for FEM and ORFD methods with different Fourier filtering levels $j^*$ and varying node counts $N$.}
\end{table}

\begin{figure}[htb!]
  \centering
	\graphicspath{{visuals/}}
	\includegraphics[width=.8\columnwidth]{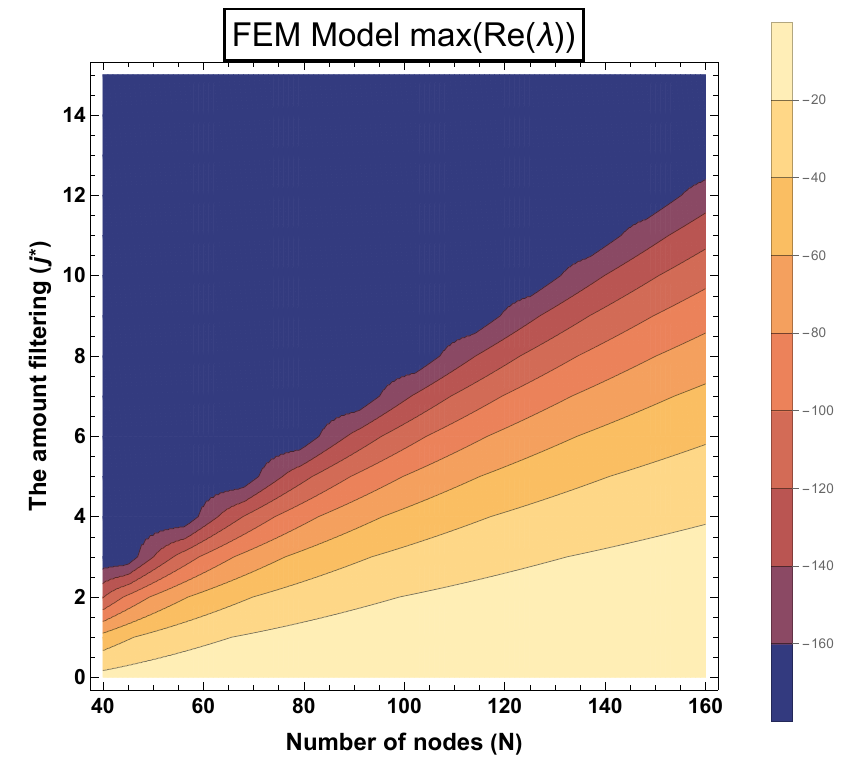}
	\caption{Maximum real part of the eigenvalues for FEM with varying Fourier filtering levels $j^*$ and different $N$ values. The total number of eigenvalues is $4N+4$, and $4j^*$ of them are filtered out.}
	\label{fig:FEMDecay}
\end{figure}

A more detailed contour plot illustrating the dependence of the maximum real part of the eigenvalues on the number of nodes and the filtering amount is presented in \cref{fig:FEMDecay}. This visualization underscores the critical role of filtering in stabilizing the system by eliminating spurious high-frequency eigenvalues with large real parts. Once a sufficient filtering threshold is reached—where lower-frequency eigenvalues dominate—further filtering no longer accelerates the decay rate. In contrast, for the ORFD case, the maximum real parts of the eigenvalues remain uniformly bounded away from the imaginary axis for all tested values of $N$. Unlike FEM solutions, ORFD solutions exhibit faster exponential decay without requiring any filtering, demonstrating their superior stability and robustness across different discretization levels.

Recalling from Section \ref{FEMSplines} that FEM solutions exhibit poor behavior for high-frequency vibrational modes, we intentionally choose high-frequency initial conditions for the simulations. For $N=80$,
\begin{equation*}
v_0 (x_j),p_0 (x_j),v_1 (x_j),p_1(x_j) = 10^{-2}\sum\limits_{k=41}^{81}jh\sin(k\pi x_j).
\end{equation*}
The final time of the simulation is set to $T_{\rm final}= 0.1$ sec

\begin{figure}[htb!]
  \centering
	\graphicspath{{visuals/}}
	\includegraphics[width=0.47\columnwidth]{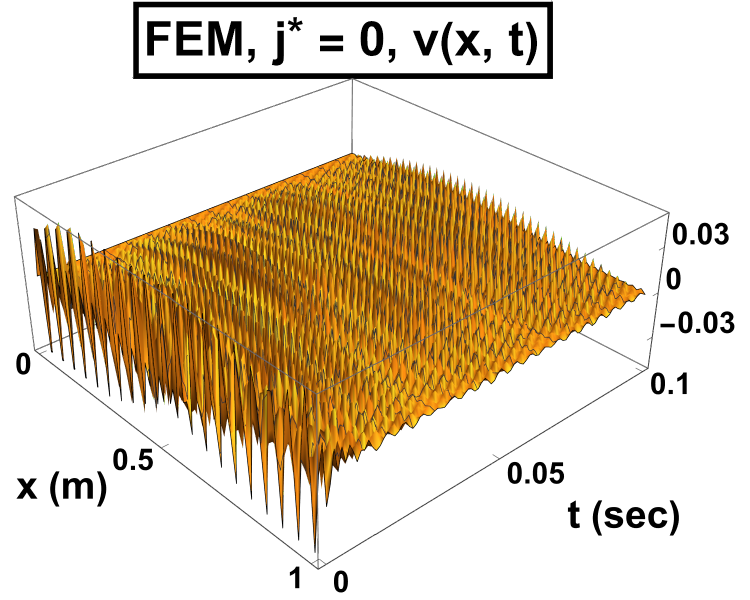} \quad \includegraphics[width=0.47\columnwidth]{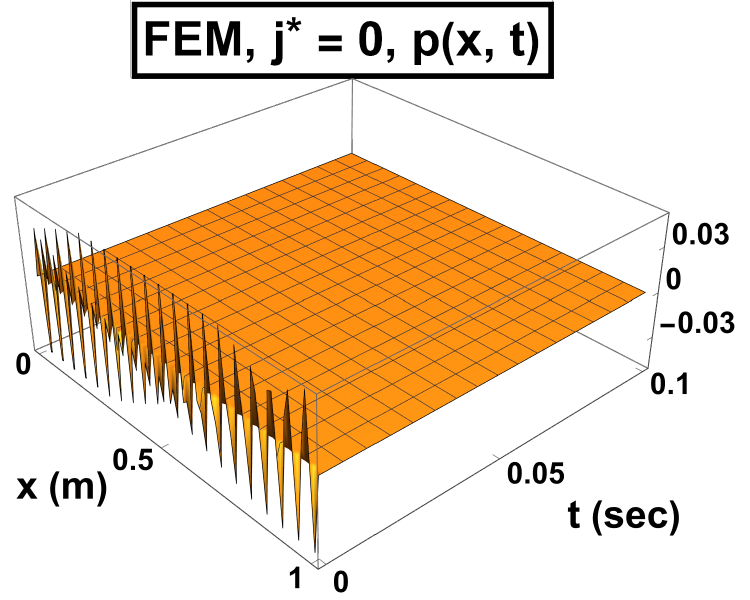}\\
	\includegraphics[width=0.47\columnwidth]{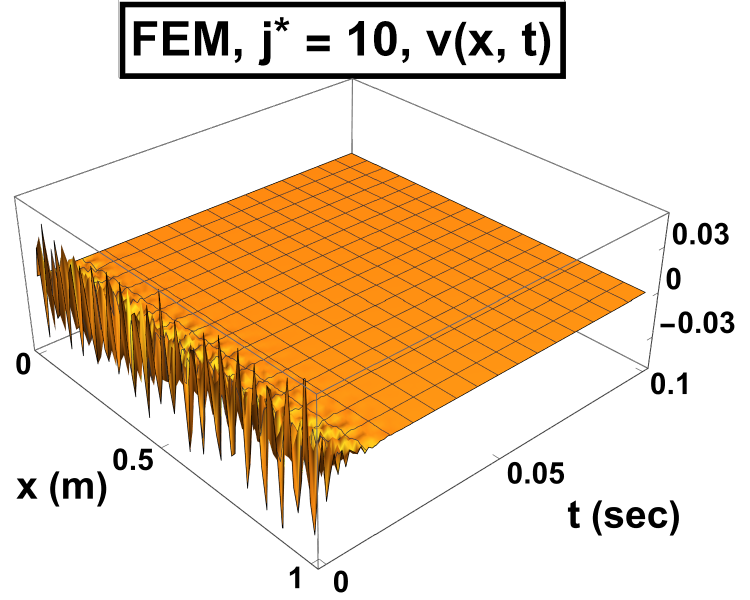} \quad \includegraphics[width=0.47\columnwidth]{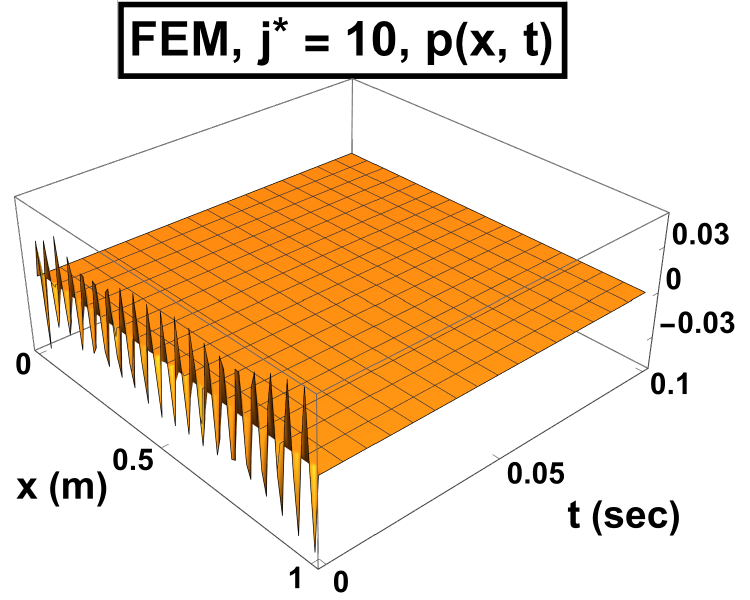}\\
	\caption{Top row: FEM solutions for $v(x,t)$ and $p(x,t)$ without filtering, showing persistent high-frequency modes. Bottom row: FEM solutions with $j^*=10$, demonstrating improved and a more realistic exponential decay behavior.}
	\label{3dFEM}
\end{figure}

\begin{figure}[htb!]
  \centering
	\graphicspath{{visuals/}}
	\includegraphics[width=0.47\columnwidth]{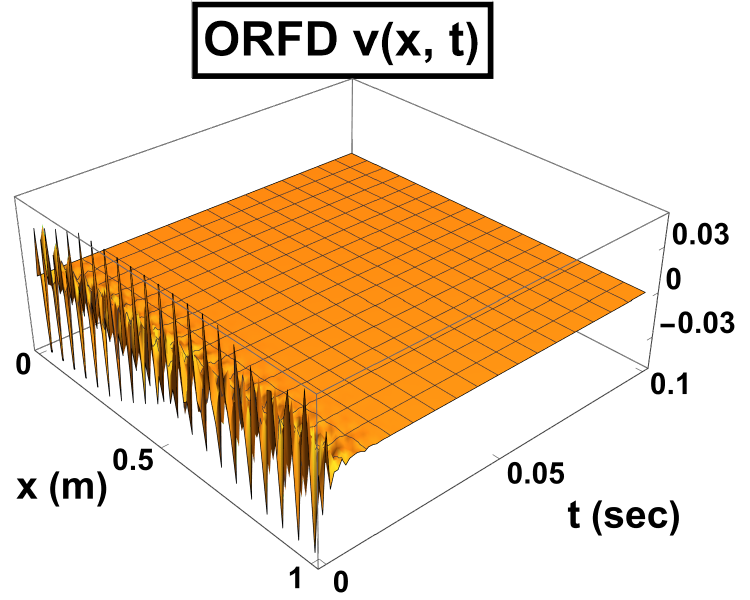} \quad \includegraphics[width=0.47\columnwidth]{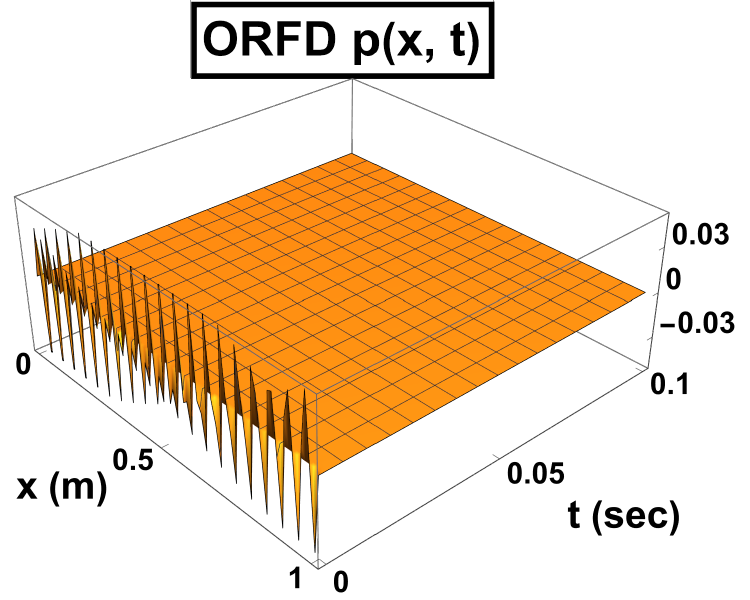}
	\caption{ORFD solutions exhibit fast exponential decay without the need for filtering, closely aligning with the PDE model.}
	\label{3dORFD}
\end{figure}

As shown in the first row of Figure \ref{3dFEM}, the FEM solutions fail to decay for high-frequency modes, leading to an unrealistic approximation of the PDE model. After applying filtering, the spurious high-frequency vibrational modes are eliminated, yielding a more accurate approximation, as depicted in the second row of Figure \ref{3dFEM}. On the other hand, the ORFD solutions exhibit exponential decay without requiring any filtering, as illustrated in Figure \ref{3dORFD}.

\begin{figure}[htb!]
  \centering
	\graphicspath{{visuals/}}
	\includegraphics[width=.8\columnwidth]{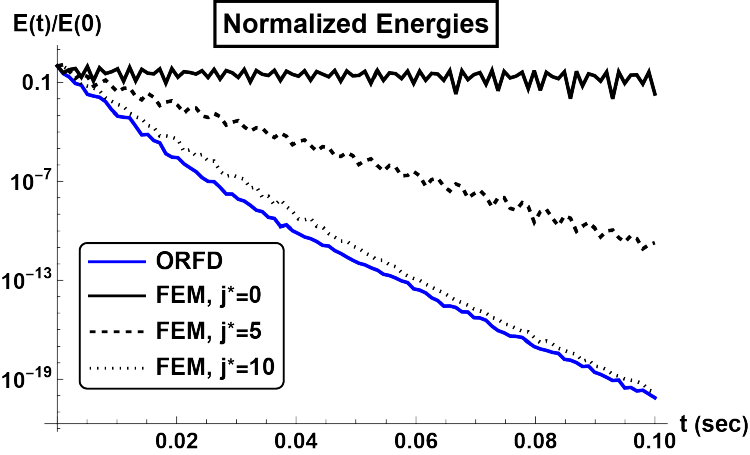}
\caption{Normalized energy decay over time for ORFD and FEM models. ORFD exhibits exponential decay without filtering, while FEM requires increasing levels of filtering (\( j^* = 5,10 \)) to achieve comparable decay rates. Without filtering (\( j^* = 0 \)), FEM solutions fail to decay properly due to persistent high-frequency modes.}
	\label{Energies}
\end{figure}

The normalized energies of the solutions are plotted on a logarithmic scale over time in Figure \ref{Energies}. These energy plots highlight the impact of high-frequency modes on stability and demonstrate the necessity of filtering in FEM models to ensure proper decay behavior.

\subsection{Influence of Feedback Amplifiers on the Optimal Filtering Parameter}

Simulations reveal the existence of a critical filtering threshold beyond which additional filtering does not enhance the decay rate. The filtering amount \( j^* \) at this threshold is referred to as the \textit{optimal filtering parameter}. A key objective is to identify this optimal \( j^* \), which ensures the fastest decay while retaining as much of the system’s dynamics as possible—i.e., applying the minimal necessary filtering to achieve the maximal decay rate.

\begin{figure}[htb!]
  \centering
	\graphicspath{{visuals/}}
	\includegraphics[width=0.8\columnwidth]{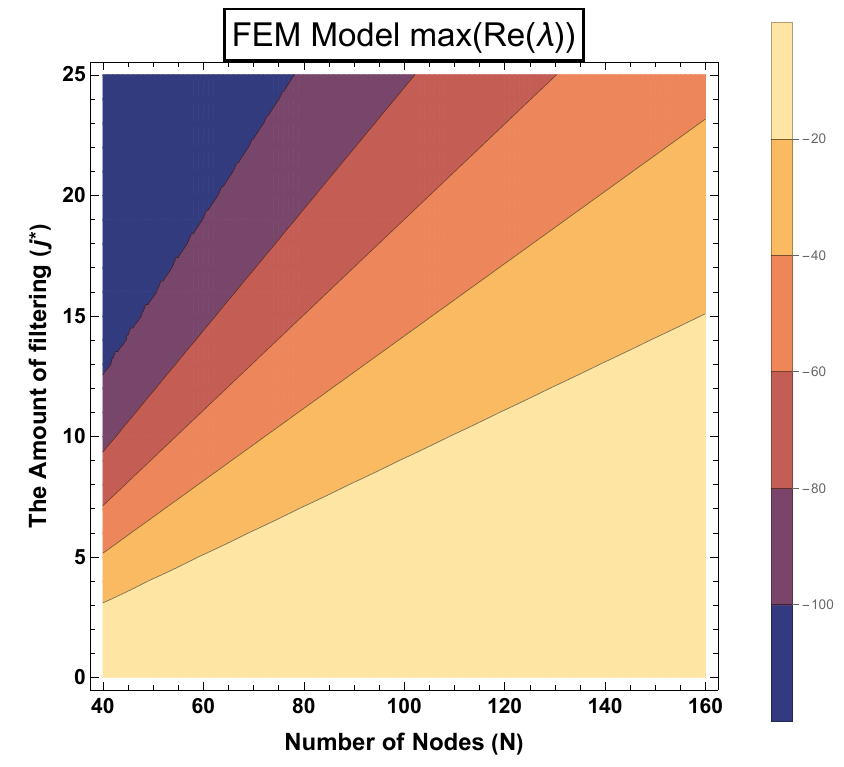}
	\caption{Maximum real part of the eigenvalues for the FEM model with varying Fourier filtering levels \( j^* \) for \( k_1 = k_2 = 10^7 \) and \( N \) nodes. The total number of eigenvalues is \( 4N+4 \), with \( 4j^* \) of them filtered out.}
	\label{fig:FEMDecayv2}
\end{figure}

Numerical results indicate a strong correlation between the feedback amplifiers and the optimal filtering parameter. Specifically, increasing the feedback gains from \( k_1 = k_2 = 10^6 \) to \( k_1 = k_2 = 10^7 \) significantly alters the filtering requirements. Figure~\ref{fig:FEMDecayv2} illustrates the dependence of the maximum real part of the eigenvalues on the Fourier filtering parameter \( j^* \) for \( N \) nodes, with \( k_1 = k_2 = 10^7 \). The total number of eigenvalues is \( 4N+4 \), with \( 4j^* \) of them filtered out.

A comparison between Figures~\ref{fig:FEMDecay} and \ref{fig:FEMDecayv2} highlights this dependency: for \( N = 80 \), the optimal filtering parameter increases from \( j^* = 5 \) (discarding approximately 6.2\% of the eigenpairs) when \( k_1 = k_2 = 10^6 \), to \( j^* = 25 \) (discarding over 31.2\% of the eigenpairs) when \( k_1 = k_2 = 10^7 \).

For further exploration, real-time Wolfram demonstrations of magnetizable piezoelectric dynamics are available at~\cite{wolframdemo}.

\section{Conclusions \& Future Work}
\label{conc}
This paper presents two novel model reduction techniques for magnetizable piezoelectric beams, addressing key challenges associated with coupled PDE systems.

Firstly, we propose a Finite Element Method (FEM)-based model reduction using linear splines, marking a new application of this approach to these systems. While this method exhibits improved stability under filtering compared to standard finite difference methods, proper application of Fourier filtering remains crucial for eliminating spurious high-frequency modes that would otherwise degrade numerical stability. However, a significant computational challenge arises: implementing Fourier filtering currently requires explicit computation of the full spectrum of the system, necessitating eigenvalue and eigenvector calculations. This process becomes computationally expensive, particularly for large-scale discretizations. An important open problem in this area is the development of efficient numerical algorithms that can filter out spurious modes without requiring explicit spectral decomposition.

Furthermore, our results reveal a fundamental connection between the optimal filtering parameter and the choice of feedback amplifiers. Numerical simulations suggest that different gain configurations necessitate varying levels of filtering to maximize stability, yet a precise theoretical characterization of this dependency remains an open question. Understanding this interplay could significantly refine the selection of control gains and filtering strategies in practical applications.

As an alternative approach, we introduce an order-reduction-based Finite Difference (ORFD) method tailored specifically to this PDE model. By leveraging Lyapunov theory, we establish exponential stability and demonstrate energy convergence of the reduced model as the discretization parameter approaches zero. Unlike the FEM-based method, ORFD does not require explicit numerical filtering, making it computationally more efficient and robust for high-fidelity simulations.

This study highlights the necessity of robust model reduction techniques, particularly in coupled systems where high-frequency eigenvalues pose computational challenges. Our approach lays the groundwork for developing model reduction strategies for multi-layer, strongly coupled systems incorporating piezoelectric layers, as explored in works such as \cite{aydin2023novel} and \cite{ozer2017modeling}-\cite{ozer2019alternate}.

Future work will extend this methodology to analyze Rao-Nakra-type three-layer beams, where wave and beam equations are strongly coupled through shear effects in the middle layer \cite{feng2023longtime}. Unlike the system considered here, which involves two branches of eigenvalues, the Rao-Nakra system exhibits three distinct eigenvalue branches, further complicating the stability and control analysis. Understanding the role of shear interactions in such systems will provide deeper insights into the interplay between different wave propagation mechanisms and contribute to more effective control and stabilization strategies for advanced smart material structures.

\bibliographystyle{abbrv}
\bibliography{references}

\begin{thebibliography}{10}

\bibitem{akil2025advancing}
M.~Akil, S.~Nicaise, A.~{\"O}. {\"O}zer, and H.~Saleh.
\newblock Advancing insights into the stabilization of novel serially-connected
  magnetizable piezoelectric and elastic beams.
\newblock {\em SIAM Journal on Control and Optimization}, 2025.
\newblock in print.

\bibitem{akil2024stability}
M.~Akil, A.~O. \"{O}zer, A.~Ramos, and H.~Wilson.
\newblock Stability results for novel serially-connected magnetizable
  piezoelectric and elastic smart-system designs.
\newblock {\em Applied Mathematics \& Optimization}, 89(3):64, 2024.

\bibitem{aydinhaider2024}
A.~K. Ayd{\i}n, M.~Z. Haider, and A.~{\"O}. {\"O}zer.
\newblock A new semi-discretization of the fully clamped euler-bernoulli beam
  preserving boundary observability uniformly.
\newblock {\em IEEE Control Systems Letters}, 2024.

\bibitem{aydin2023novel}
A.~K. Ayd{\i}n, A.~{\"O}. {\"O}zer, and J.~Walterman.
\newblock A novel finite difference-based model reduction and a sensor design
  for a multilayer smart beam with arbitrary number of layers.
\newblock {\em IEEE Control Systems Letters}, 7:1548--1553, 2023.

\bibitem{banks2012functional}
H.~T. Banks.
\newblock {\em A functional analysis framework for modeling, estimation and
  control in science and engineering}.
\newblock CRC Press, 2012.

\bibitem{banks1991exponentially}
H.~T. Banks, K.~Ito, and C.~Wang.
\newblock Exponentially stable approximations of weakly damped wave equations.
\newblock In {\em Estimation and Control of Distributed Parameter Systems:
  Proceedings of an International Conference on Control and Estimation of
  Distributed Parameter Systems, Vorau, July 8--14, 1990}, pages 1--33.
  Springer, 1991.

\bibitem{baur2014advances}
C.~Baur, D.~J. Apo, D.~Maurya, S.~Priya, and W.~Voit.
\newblock Advances in piezoelectric polymer composites for vibrational energy
  harvesting.
\newblock In {\em Polymer composites for energy harvesting, conversion, and
  storage}, pages 1--27. ACS Publications, 2014.

\bibitem{darinskii2008role}
A.~Darinskii, E.~Le~Clezio, and G.~Feuillard.
\newblock The role of electromagnetic waves in the reflection of acoustic waves
  in piezoelectric crystals.
\newblock {\em Wave Motion}, 45(4):428--444, 2008.

\bibitem{de2023control}
M.~C. de~Jong, K.~C. Kosaraju, and J.~M. Scherpen.
\newblock On control of voltage-actuated piezoelectric beam: A krasovskii
  passivity-based approach.
\newblock {\em European Journal of Control}, 69:100724, 2023.

\bibitem{el2013boundary}
H.~El~Boujaoui, H.~Bouslous, and L.~Maniar.
\newblock Boundary stabilization for 1-d semi-discrete wave equation by
  filtering technique.
\newblock {\em Bull. TICMI}, 17(1):1--18, 2013.

\bibitem{feng2022exponential}
B.~Feng and A.~{\"O}. {\"O}zer.
\newblock Exponential stability results for the boundary-controlled
  fully-dynamic piezoelectric beams with various distributed and boundary
  delays.
\newblock {\em Journal of Mathematical Analysis and Applications},
  508(1):125845, 2022.

\bibitem{feng2023longtime}
B.~Feng and A.~{\"O}. {\"O}zer.
\newblock Long-time behavior of a nonlinearly damped rao-nakra sandwich beam.
\newblock {\em Applied Mathematics and Optimization}, 87:Article number: 19 (52
  pages), 2023.

\bibitem{feng2023stability}
B.~Feng and A.~{\"O}. {\"O}zer.
\newblock Stability results for piezoelectric beams with long-range memory
  effects in the boundary.
\newblock {\em Mathematische Nachrichten}, 296(9):4206--4235, 2023.

\bibitem{freitas2022long}
M.~Freitas, A.~{\"O}. {\"O}zer, and A.~Ramos.
\newblock Long time dynamics and upper semi-continuity of attractors for
  piezoelectric beams with nonlinear boundary feedback.
\newblock {\em ESAIM: Control, Optimisation and Calculus of Variations}, 28:39,
  2022.

\bibitem{infante1999boundary}
J.~A. Infante and E.~Zuazua.
\newblock Boundary observability for the space semi-discretizations of the 1--d
  wave equation.
\newblock {\em ESAIM: Mathematical Modelling and Numerical Analysis},
  33(2):407--438, 1999.

\bibitem{kong2022equivalence}
A.~Kong, C.~Nonato, W.~Liu, M.~J.~d. Santos, and C.~Raposo.
\newblock Equivalence between exponential stabilization and observability
  inequality for magnetic effected piezoelectric beams with time-varying delay
  and time-dependent weights.
\newblock {\em Discrete \& Continuous Dynamical Systems-Series B}, 27(6), 2022.

\bibitem{leon2002boundary}
L.~Le{\'o}n and E.~Zuazua.
\newblock Boundary controllability of the finite-difference space
  semi-discretizations of the beam equation.
\newblock {\em ESAIM: Control, Optimisation and Calculus of Variations},
  8:827--862, 2002.

\bibitem{lissy2019optimal}
P.~Lissy and I.~Roven{\c{t}}a.
\newblock Optimal filtration for the approximation of boundary controls for the
  one-dimensional wave equation using a finite-difference method.
\newblock {\em Mathematics of Computation}, 88(315):273--291, 2019.

\bibitem{liu2019novel}
J.~Liu and B.-Z. Guo.
\newblock A novel semi-discrete scheme preserving uniformly exponential
  stability for an euler--bernoulli beam.
\newblock {\em Systems \& Control Letters}, 134:104518, 2019.

\bibitem{liu2020new}
J.~Liu and B.-Z. Guo.
\newblock A new semidiscretized order reduction finite difference scheme for
  uniform approximation of one-dimensional wave equation.
\newblock {\em SIAM Journal on Control and Optimization}, 58(4):2256--2287,
  2020.

\bibitem{liu2021uniformly}
J.~Liu and B.-Z. Guo.
\newblock Uniformly semidiscretized approximation for exact observability and
  controllability of one-dimensional euler--bernoulli beam.
\newblock {\em Systems \& Control Letters}, 156:105013, 2021.

\bibitem{morris2014modeling}
K.~A. Morris and A.~{\"O}. {\"O}zer.
\newblock Modeling and stabilizability of voltage-actuated piezoelectric beams
  with magnetic effects.
\newblock {\em SIAM Journal on Control and Optimization}, 52(4):2371--2398,
  2014.

\bibitem{ozer2015further}
A.~{\"O}. {\"O}zer.
\newblock Further stabilization and exact observability results for
  voltage-actuated piezoelectric beams with magnetic effects.
\newblock {\em Mathematics of Control, Signals, and Systems}, 27(2):219--244,
  2015.

\bibitem{ozer2017modeling}
A.~{\"O}. {\"O}zer.
\newblock Modeling and controlling an active constrained layer (acl) beam
  actuated by two voltage sources with/without magnetic effects.
\newblock {\em IEEE Transactions on Automatic Control}, 62(12):6445--6450,
  2017.

\bibitem{ozer2018potential}
A.~{\"O}. {\"O}zer.
\newblock Potential formulation for charge or current-controlled piezoelectric
  smart composites and stabilization results: electrostatic versus quasi-static
  versus fully-dynamic approaches.
\newblock {\em IEEE Transactions on Automatic Control}, 64(3):989--1002, 2018.

\bibitem{ozer2019uniform}
A.~{\"O}. {\"O}zer.
\newblock Uniform boundary observability of semi-discrete finite difference
  approximations of a rayleigh beam equation with only one boundary
  observation.
\newblock In {\em 2019 IEEE 58th Conference on Decision and Control (CDC)},
  pages 7708--7713. IEEE, 2019.

\bibitem{ozer2021stabilization}
A.~{\"O}. {\"O}zer.
\newblock Stabilization results for well-posed potential formulations of a
  current-controlled piezoelectric beam and their approximations.
\newblock {\em Applied Mathematics \& Optimization}, 84(1):877--914, 2021.

\bibitem{ozer2022robust}
A.~{\"O}. {\"O}zer and A.~K. Ayd{\i}n.
\newblock Robust-filtering of sensor data for the finite difference model
  reduction of a piezoelectric sandwich beam.
\newblock In {\em 2022 IEEE 61st Conference on Decision and Control (CDC)},
  pages 6535--6541. IEEE, 2022.

\bibitem{ozer2024exponential}
A.~{\"O}. {\"O}zer, A.~K. Ayd{\i}n, and R.~Emran.
\newblock Exponential stability and design of sensor feedback amplifiers for
  fast stabilization of magnetizable piezoelectric beam equations.
\newblock {\em IEEE Transactions on Automatic Control}, 2024.

\bibitem{ozer2023revisiting}
A.~{\"O}. {\"O}zer and R.~Emran.
\newblock Revisiting the direct fourier filtering technique for the maximal
  decay rate of boundary-damped wave equation by finite differences and finite
  elements.
\newblock {\em under revision, available at arXiv:2306.11398}, 2023.

\bibitem{ozer2022uniform}
A.~{\"O}. {\"O}zer and W.~Horner.
\newblock Uniform boundary observability of finite difference approximations of
  non-compactly coupled piezoelectric beam equations.
\newblock {\em Applicable Analysis}, 101(5):1571--1592, 2022.

\bibitem{ozer2019alternate}
A.~{\"O}. {\"O}zer and M.~Khenner.
\newblock An alternate numerical treatment for nonlinear pde models of
  piezoelectric laminates.
\newblock In {\em Active and Passive Smart Structures and Integrated Systems
  XIII}, volume 10967, pages 441--460. SPIE, 2019.

\bibitem{ozer2024boundary}
A.~{\"O}. {\"O}zer, U.~Rasaq, and I.~Khalilullah.
\newblock Boundary output feedback stabilization for a novel magnetizable
  piezoelectric beam model.
\newblock In {\em IEEE American Control Conference (ACC)}, pages 3448--3453.
  IEEE, 2024.

\bibitem{preumont2018vibration}
A.~Preumont.
\newblock {\em Vibration control of active structures: an introduction}, volume
  246.
\newblock Springer, 2018.

\bibitem{ramos2019equivalence}
A.~Ramos, M.~Freitas, D.~Almeida, S.~Jesus, and T.~Moura.
\newblock Equivalence between exponential stabilization and boundary
  observability for piezoelectric beams with magnetic effect.
\newblock {\em Zeitschrift f{\"u}r angewandte Mathematik und Physik}, 70:1--14,
  2019.

\bibitem{ramos2018exponential}
A.~J. Ramos, C.~S. Gon{\c{c}}alves, and S.~S. Corr{\^e}a~Neto.
\newblock Exponential stability and numerical treatment for piezoelectric beams
  with magnetic effect.
\newblock {\em ESAIM: Mathematical Modelling and Numerical Analysis},
  52(1):255--274, 2018.

\bibitem{smith2005smart}
R.~C. Smith.
\newblock {\em Smart material systems: model development}.
\newblock SIAM, 2005.

\bibitem{tallarico2017propagation}
D.~Tallarico, N.~Movchan, A.~Movchan, and M.~Camposaragna.
\newblock Propagation and filtering of elastic and electromagnetic waves in
  piezoelectric composite structures.
\newblock {\em Mathematical Methods in the Applied Sciences}, 40(9):3202--3220,
  2017.

\bibitem{tebou2007uniform}
L.~T. Tebou and E.~Zuazua.
\newblock Uniform boundary stabilization of the finite difference space
  discretization of the 1- d wave equation.
\newblock {\em Advances in Computational Mathematics}, 26:337--365, 2007.

\bibitem{wolframdemo}
J.~Walterman, A.~K. Ayd{\i}n, S.~Leveridge, and A.~{\"O}. {\"O}zer.
\newblock Dynamics of a longitudinal piezoelectric beam.
\newblock {\em Wolfram Demonstrations Project}, 2023.

\bibitem{yang2009fully}
J.~Yang.
\newblock Fully dynamic theory.
\newblock In {\em Special Topics in the Theory of Piezoelectricity}, pages
  247--280. Springer, 2009.

\bibitem{zhang2023stability}
H.-E. Zhang, G.-Q. Xu, and Z.-J. Han.
\newblock Stability and eigenvalue asymptotics of multi-dimensional fully
  magnetic effected piezoelectric system with friction-type infinite memory.
\newblock {\em SIAM Journal on Applied Mathematics}, 83(2):510--529, 2023.

\bibitem{zuazua2006controllability}
E.~Zuazua.
\newblock {\em Controllability of partial differential equations}.
\newblock PhD thesis, Optimization and Control, 2006.

\end{thebibliography}

\medskip

\end{document}